\documentclass[final,leqno,onefignum,onetabnum]{siamltex}
\pdfoutput=1
\usepackage{amsmath}
\usepackage{xcolor}
\usepackage{amssymb}
\usepackage{amsfonts}
\usepackage{fancybox}
\usepackage{color}
\usepackage{bm}
\usepackage{cancel}
\usepackage{framed}

\newtheorem{remark}{Remark}

\newcommand{\mauro}[1]{{\color{blue}{#1}}}

\usepackage{graphicx}
\usepackage{epstopdf}
\DeclareMathOperator*{\argmin}{argmin}
\definecolor{shadecolor}{gray}{0.75}
\newcommand{\email}[1]{\protect\href{mailto:#1}{#1}}

\title{A high-order staggered meshless method for elliptic problems}

\author{Nathaniel Trask
\thanks{Division of Applied Mathematics, Brown University, 182 George St. Providence, RI 02906 (\email{nat.trask@gmail.com})}
\and
Mauro Perego\thanks{Center for Computing Research, Sandia National Laboratories, Mail Stop 1320 Albuquerque, New Mexico, 87185-1320 ({\tt \{pbboche,mperego\}@sandia.gov}). Sandia National Laboratories is a multi-program laboratory managed and operated by Sandia Corporation, a wholly owned subsidiary of Lockheed Martin Corporation, for the U.S. Department of Energy's National Nuclear Security Administration under contract DE-AC04-94AL85000} 
\and 
Pavel Bochev$^\dagger$}

\begin{document}
\maketitle
\slugger{sisc}{xxxx}{xx}{x}{x--x}

\begin{abstract}
We present a new meshless method for scalar diffusion equations, which is motivated by their compatible discretizations on primal-dual grids. 
Unlike the latter though, our  approach is truly meshless because it only requires the graph of nearby neighbor connectivity of the discretization points $\bm{x}_i$. This graph defines a local primal-dual grid complex with a \emph{virtual} dual grid, in the sense that specification of the dual metric attributes is implicit in the method's construction. 
Our method combines a topological gradient operator on the local primal grid with a Generalized Moving Least Squares approximation of the divergence on the local dual grid. 
We show that the resulting approximation of the div-grad operator maintains polynomial reproduction to arbitrary orders and  yields a meshless method, which attains $O(h^{m})$ convergence in both $L^2$ and $H^1$ norms, similar to mixed finite element methods.
We demonstrate this convergence on curvilinear domains using manufactured solutions.
Application of the new method to problems with discontinuous coefficients reveals solutions that are qualitatively similar to those of compatible mesh-based discretizations.
\end{abstract}


\begin{AMS}\end{AMS}

\pagestyle{myheadings}
\thispagestyle{plain}
\markboth{N. TRASK et al.}{Staggered moving least squares}

\section{Introduction}
Consider a bounded region in $\mathbb{R}^d$, $d=2,3$ with a Lipschitz continuous boundary $\Gamma=\partial\Omega$. This paper presents a new staggered meshless discretization approach for the model elliptic boundary value problem
\begin{equation}\label{eq:model}
\begin{array}{rcl}
-\nabla\cdot\mu\nabla \phi & = & f \quad \mbox{in $\Omega$} \\[1ex]
\phi &=&  u\quad \mbox{on $\Gamma_D$} \\[1ex]
\bm{n}\cdot\mu\nabla\phi &=& g \quad \mbox{on $\Gamma_N$} 
\end{array}
\end{equation}
where $\Gamma_D$ and $\Gamma_N$ denote Dirichlet and Neumann parts of the boundary $\Gamma$, respectively, $\mu$ is a symmetric positive definite tensor describing a material property and $f$, $u$ and $g$ are given data.

In important applications such as porous media flow \cite{yotov-phd}, heat transfer \cite{Swaminathan_93_FEAD} and semiconductor devices \cite{Brezzi_89_CMAME},  the flux $\bm{u}=-\mu\nabla \phi$ is the variable of primary interest. Predictive simulations of such problems require carefully constructed approximations of the divergence and gradient operators, which ensure stable, accurate and locally conservative flux approximations. 

Generally speaking, such \emph{compatible discretizations} of \eqref{eq:model}
fall into one of the following two categories. \emph{Single grid} methods such as: mixed finite elements \cite{Brezzi_91_BOOK}; virtual elements \cite{DaVeiga_13_MMMAS}; or mimetic finite differences  \cite{Lipnikov_06_JCP}, \cite{Shashkov_95_BOOK},  
construct an approximation of one of the two operators (typically the divergence) on the given grid and then use its adjoint\footnote{Multiplication of $\bm{u}=-\mu\nabla \phi$ by a test function $\bm{v}$ followed by integration by parts yields a notion of a weak gradient in mixed finite elements. Other methods, such as mimetic finite differences define the gradient as the algebraic adjoint of the discrete divergence.} as a proxy for the other operator (the gradient). In so doing these methods satisfy a discrete version of Green's theorem. 

In contrast, \emph{primal-dual grid} approaches, such as box integration \cite{Kramer_97_BOOK}, covolume  \cite{Trapp_06_THESIS,Nicolaides_97_SINUM}, or finite volume \cite{Mishev_07_NMPDE} schemes  use the two grids in the primal-dual complex to approximate the divergence and the gradient independently. When the two grids are topologically dual, such as in the case of Voronoi-Delauney triangulations or rectilinear\footnote{In this case the dual grid is a translation of the primal grid by a half cell size. This creates the appearance of all variables living on the same grid but at different locations, thus the term ``staggered grid methods''.} primal grids, these methods assume a particularly simple and elegant form.

The abundance of compatible mesh-based methods that provide accurate and locally conservative flux approximations stands in stark contrast with the almost complete absence of meshless methods that would provide comparable results but without the need for a global mesh structure. 
Few notable exceptions are the uncertain grid method \cite{Diyankov_08_TECHREP} and the conceptually similar meshless volume schemes \cite{Katz_09_AIAA,chiu2011conservative} for conservation laws. 

Motivated by compatible discretizations of \eqref{eq:model} on primal-dual grids, we present a new meshless scheme for this problem that exhibits similar computational properties, most notably an ability to produce accurate, non-oscillatory approximations for problems with material discontinuities.
Our  approach is truly meshless because it only requires the graph of nearby neighbor connectivity of every discretization point $\bm{x}_i$. The vertices and the edges of this graph define a local primal grid for every mesh point, which induces a virtual local dual grid comprising a virtual cell dual to $\bm{x}_i$ and a set of virtual faces dual to the edges of the connectivity graph. 

Following primal-dual grid methods we discretize the div-grad operator using independent approximations of the divergence and the gradient on the local grid complex. Specifically, we combine a topological gradient on the local primal grid with a Generalized Moving Least Squares (GMLS) approximation of the divergence on the local virtual dual grid. 
The latter is virtual in the sense that our approach does not require its physical construction, instead, specification of dual cell volumes and dual face areas is implicit in the method's construction. 
We show that the resulting meshless approximation of the div-grad operator retains polynomial reproduction of arbitrary orders. Because of its conceptual similarity with primal-dual schemes  we call the resulting method a \emph{staggered GMLS} discretization of \eqref{eq:model}.

Our new method also differs in important and significant ways from the conservative meshless schemes in \cite{Diyankov_08_TECHREP,Katz_09_AIAA,chiu2011conservative}. While these approaches similarly use 
different sets of degrees-of-freedom to discretize the gradient and divergence, they rely on a global Quadratic Program (QP) to determine the meshless analogues of cell volumes and face areas that ensure the compatibility of the scheme. As a result, application of these methods to large-scale problems would require the formulation of a scalable, QP-specific optimization solver \cite{chiu2011conservative}. In contrast, our approach requires the solution of many independent and inexpensive local optimization problems, allowing an $O(N)$ implementation in the number of degrees of freedom.

A second important difference is the order of accuracy achieved by these schemes and our new method. The accuracy of the former is similar to that of comparable finite volume schemes, i.e., they are roughly second-order accurate. In contrast, our method attains $O(h^{m})$ convergence in both $L^2$ and $H^1$ norms, similar to mixed finite element methods.

It is worth mentioning that some variations of smooth particle hydrodynamics (SPH), which use discretely skew-adjoint divergence and gradient operators \cite{monaghan2012smoothed,price2008modelling}, pursue a somewhat different type of compatibility. However, it is well-known \cite{traskCMAME} that this compatibility in SPH comes at the expense of the gradient operators lacking even zeroth-order polynomial consistency, and many applications require techniques to numerically stabilize the method in the presence of large differences in material properties such as density or viscosity (e.g. \cite{hu2007incompressible}). In contrast to this, we will show that our method maintains high-order polynomial reproduction and remains numerically stable across jumps in material properties of several orders of magnitude.

We have organized the paper as follows. Section \ref{sec:not} introduces notation and reviews basic concepts of Generalized Moving Least Squares, while Section \ref{sec:GMLS-div} specializes these concepts to approximation of vector fields and their derivatives from directional components. 
Section \ref{sec:define} presents the development of the new staggered GMLS, and Section \ref{sec:props} studies its accuracy. We discuss implementation details in Section \ref{sec:implement} and present numerical results in Section \ref{sec:results}. Section \ref{sec:conclude} summarizes our findings.

\section{Notation and quotation of results}\label{sec:not}
Throughout this paper bold face fonts denote various vector quantities, e.g., $\bm{x}$ is a point in the Euclidean space $\mathbb{R}^d$, $\bm{u}$ is a vector field in $\mathbb{R}^d$, $\bm{e}_{ij}$ is an edge connecting points $\bm{x}_i$ and $\bm{x}_j$, and so on. Lower case fonts stand for various scalar quantities, upper case Roman fonts are primarily used to denote function spaces, operators and sets of geometric entities, while the first few lower case Greek symbols will stand for multi-indices, e.g., $\alpha=(\alpha_1,\ldots,\alpha_d)$ and $|\alpha | = \sum\alpha_i$.
We use the standard Euler notation $D^\alpha$ for a partial derivative of order $|\alpha|$, i.e.,
 $D^\alpha:=\partial^{|\alpha|} /\partial_{\bm{x}_{i}^{\alpha_1}}\ldots\partial_{\bm{x}_{i}^{\alpha_i}}$ and denote the standard Euclidean norm by $\|\cdot\|$. 

As usual, $C^m(\Omega)$ is the space of all continuous functions whose derivatives up to order $m$ are also continuous and $P_m(\mathbb{R}^d)$, or simply $P_m$ is the space of all multivariate polynomials of degree less than or equal to $m$.

\subsection{Local geometric structures}\label{sec:local-geom}

\begin{figure}[h!]
  \centering
   \includegraphics[height=2in]{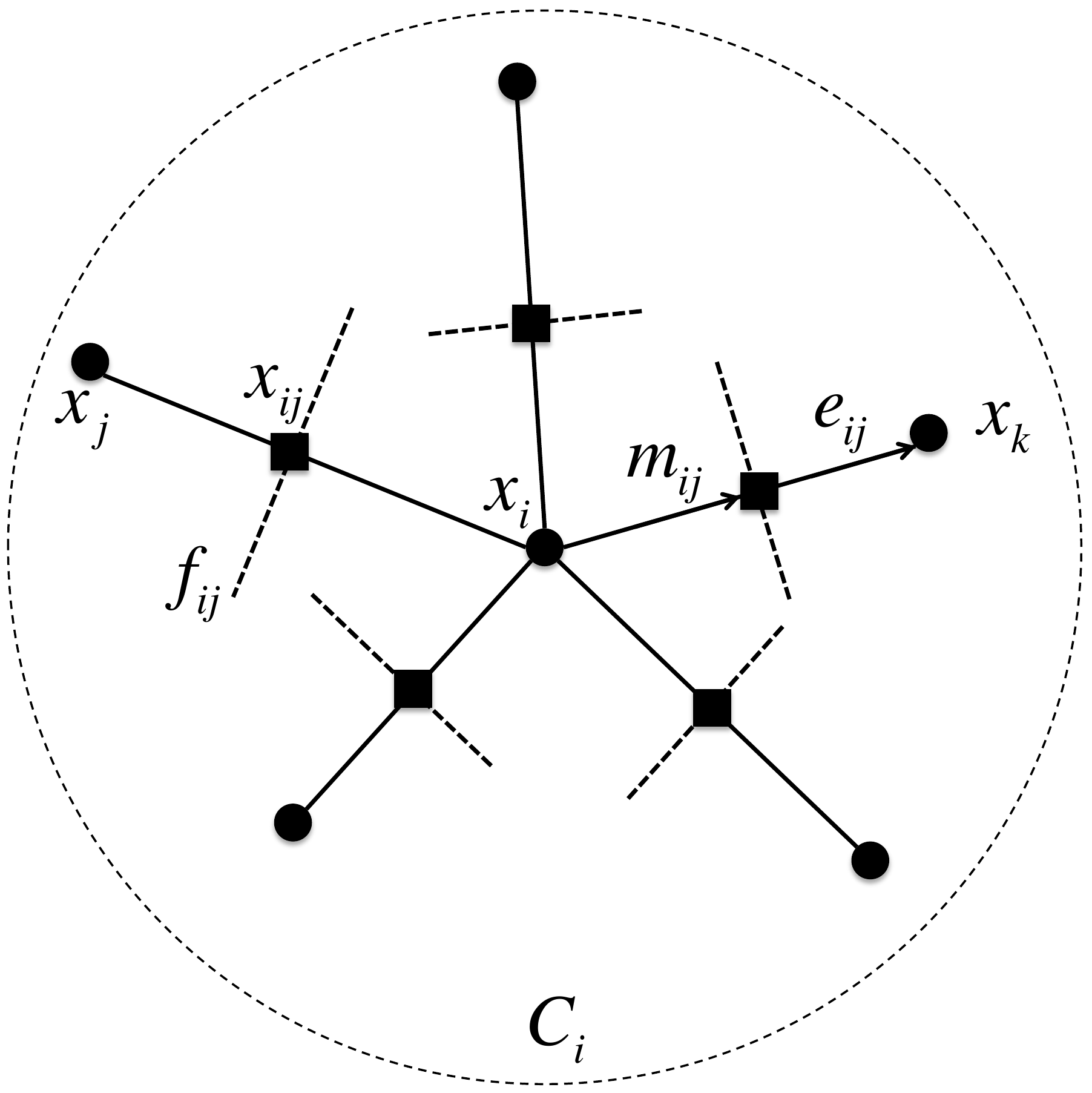} 
  \caption{A local primal-dual grid complex induced by a point $\bm{x}_i$ and its $\varepsilon$-neighborhood $N^\varepsilon_i$. The primal edges $\bm{e}_{ij}$, midpoints $\bm{x}_{ij}$ and mid-edges $\bm{m}_{ij}$ are physical mesh entities. The dual cell $C_i$ and the dual faces $\bm{f}_{ij}$ are virtual mesh entities.}
   \label{stencil}
\end{figure}

In this paper we consider discretization of $\Omega$ by a set of points $\Omega_h=\{\bm{x}_i\}_{i=1,\ldots,N}$. We denote all boundary points by $\Gamma_h$ and set 
$ \breve{\Omega}_h = \Omega_h\setminus \Gamma_h$.
The $\varepsilon$-neighborhood of $\bm{x}_i\in\Omega_h$ is the set
$$
N^\varepsilon_i:=\{\bm{x}_j\in\Omega_h\, |\, \|\bm{x}_j - \bm{x}_i\|<\varepsilon\},
$$
where $\varepsilon>0$ is given and may depend on the point location. 
We use the local connectivity graph of the points in $N^\varepsilon_i$ to  construct an approximation of the div-grad operator that mimics the one in primal-dual grid methods. 
Specifically, with every $\varepsilon$-neighborhood we associate a \emph{local primal grid} comprising a local vertex set
$$
V_i = \{\bm{v}_j=\bm{x}_j\, |\, \bm{x}_j\in N^\varepsilon_i\},
$$
which is simply the set of all points in $N^\varepsilon_i$, and a local edge set
$$
E_i = \{\bm{e}_{ij}=\bm{x}_{j}-\bm{x}_i \,|\,\bm{x}_j\in N^\varepsilon_i\}.
$$
The edges in $E_i$ have midpoints, half-edges, and unit tangents given by
$$
\bm{x}_{ij} = \frac{\bm{x}_i+\bm{x}_j}{2},\ \
\bm{m}_{ij} =  \bm{x}_{ij}-\bm{x}_i,
\quad\mbox{and}\quad
\bm{t}_{ij} =  \frac{\bm{x}_j-\bm{x}_i}{\|\bm{x}_j-\bm{x}_i\|} = \frac{\bm{x}_j-\bm{x}_i}{\|\bm{e}_{ij}\|},
$$
respectively; see Fig.~\ref{stencil}. We denote the set of all midpoints by $M_i$. Note that $E_i$ does not contain edges that do not have $\bm{x}_i$ as a vertex.

The local primal grid induces a \emph{local virtual dual grid} comprising a single virtual cell $C_i$ dual to vertex $\bm{v}_i$ and a set of virtual faces
$
F_i = \{\bm{f}_{ij} \}
$
dual to $E_i$, intersecting $\bm{e}_{ij}$ at $\bm{x}_{ij}$ and having face normals equivalent to edge tangents, i.e.,
\begin{equation}\label{eq:topodual}
\forall \bm{f}_{ij} \in F_i,\quad \bm{f}_{ij} \cap \bm{e}_{ij} = \bm{x}_{ij} \quad \mbox{and} \quad \bm{n}_{ij} = \bm{t}_{ij}\,.
\end{equation}

\begin{remark}
Although the local primal mesh contains ``real'' mesh entities in the sense that the elements of $V_i$ and $E_i$ have geometric attributes such as coordinates and lengths, they are not used directly in the formulation of our method. Instead, all necessary metric quantities are determined implicitly by the method's construction.
\end{remark}

\subsection{Local approximation spaces and operators}
On the local primal grid we define the local vertex space
$$
V^i=\{ u_j\in\mathbb{R}\,|\, \forall \bm{v}_j\in V_i\}
$$ 
and the local edge space
$$
E^i=\{ u_{ij}\in\mathbb{R}\,|\, \forall \bm{e}_{ij}\in E_i\}\,.
$$ 
The elements of $V^i$ and $E^i$ are sets of real values associated with the local vertices and edges, respectively. In particular, we assume that edge values $u_{ij}$ are localized at the edge midpoints $\bm{x}_{ij}$ and represent the tangent components of some vector field $\bm{u}$ along the edges, i.e., $u_{ij} = \bm{u}_{ij}\cdot\bm{t}_{ij}$.

Likewise, on the local virtual dual grid we define the cell space
$$
C^i = \{ u_i\in\mathbb{R}\}
$$
comprising a single scalar value on the virtual dual cell $C_i$, and the face space
$$
F^i = \{ u_{ij}\in\mathbb{R}\,|\, \forall \bm{f}_{ij}\in F_i\}
$$
containing real values associated with the virtual dual faces. We assume that the face values $u_{ij}$ are localized at edge-face intersection points $\bm{x}_{ij}$ and represent the normal components of some vector field $\bm{u}$, i.e., $u_{ij} = \bm{u}_{ij}\cdot\bm{n}_{ij}$.

Thanks to the topological duality  \eqref{eq:topodual} between $E_i$ and $F_i$, the corresponding edge and face spaces $E^i$ and $F^i$ are isomorphic. Thus, any tangent vector component $u_{ij}=\bm{u}_{ij}\cdot\bm{t}_{ij}$ on a primal edge $\bm{e}_{ij}$ can be viewed as a normal vector component $\bm{u}_{ij}\cdot\bm{n}_{ij}$ on the dual face $\bm{f}_{ij}$ and vice versa.  Thus, our local approximation spaces reproduce the key property of primal-dual grid complexes without requiring a global mesh data structure. 

To complete the specification of the discrete structures necessary for our new scheme it remains to endow the above discrete spaces with suitable notions of a gradient and divergence operators. Taking clue from primal-dual grid methods we seek the discrete gradient as a mapping  $GRAD_{i} : V^i\rightarrow E^i$ on the local primal mesh, and the discrete divergence as a mapping $DIV_i : F^i\rightarrow C^i$ on the local virtual dual mesh. 
As in a primal-dual  methods, the isomorphism between $E^i$ and $F^i$ then allows us to ``chain'' these operators into an approximation of the div-grad operator. 

Avoiding an explicit dependence on any metric entities such as areas, volumes and lengths is a key requirement for the construction of $GRAD_{i}$ and $DIV_i$. This enables us to bypass a global optimization problem as required in \cite{chiu2011conservative}. To this end, we choose to define $GRAD_i$ as the \emph{topological} gradient \cite{Bochev_06a_IMA,Hyman_97_CMA}, i.e., we set
\begin{equation}\label{eq:grad}
GRAD_i(u^h) = \bm{u}^h \in E^i 
\Leftrightarrow
u_{ij} = u_j - u_i
\quad \forall u^h\in V^i .
\end{equation}
The operator \eqref{eq:grad} is topological because it depends on the connectivity of the local vertices and edges in $N^\varepsilon_i$, but not on their physical locations. The node-to-edge incidence matrix of the local primal grid is the algebraic representation of $GRAD_i$. 

We will subsume all necessary metric attributes of our scheme in the definition of  $DIV_i$ by using a Generalized Moving Least Squares approach \cite{mirzaei2011} to define this operator. We now review key aspects of the GMLS theory necessary for this task. Then in Section \ref{sec:GMLS-div} we extend the approach in \cite{mirzaei2011} to approximate vector fields and their derivatives from scattered directional components

\subsection{Generalized Moving Least Squares (GMLS) framework}\label{sec:GMLS}
Following \cite[Section 4.3]{WendlandBook} we consider an abstract setting given by
\begin{itemize}
\item a function space $V$ with a dual $V^*$;
\item a finite dimensional space $P=\mbox{span}\{p_1,\ldots,p_Q\}\subset V$; 
\item a finite set of linear functionals $\Lambda=\{\lambda_1,\ldots,\lambda_N\}\subset V^*$; and 
\item a correlation (weight) function $\omega:V^*\times V^*\mapsto \mathbb{R}^+\cup\{0\}$.
\end{itemize}
We assume that $\Lambda$ is $P$-unisolvent, that is
\begin{equation}\label{eq:P-uni}
\{p\in P\,|\, \lambda_i(p) = 0, i=1,\ldots,N\} = \{0\}.
\end{equation}
In other words, the zero is the only element of $P$ for which all functionals in $\Lambda$ vanish.

Given a \emph{target} functional $\tau\in V^*$ GMLS seeks to approximate its action on any function $u\in V$ by a linear combination 
\begin{equation}\label{eq:gmls}
\tau(u) \approx \widetilde{\tau}(u):= \sum_{i=1}^{N} a_i(\tau) \lambda_i(u)
\end{equation}
of its samples $\{\lambda_1(u),\ldots,\lambda_N(u)\}$ such that 

\noindent
\textbf{P.1} The approximation \eqref{eq:gmls} is exact for $P$, i.e.,
\begin{equation}\label{eq:gmls-ex}
\tau(u) = \widetilde{\tau}(u):= \sum_{i=1}^{N} a_i(\tau) \lambda_i(u)\quad\forall u\in P.
\end{equation}

\noindent
\textbf{P.2} The coeffcients $a_i(\tau)$ have local supports relative to $\omega$, i.e.,
\begin{equation}\label{eq:gmls-loc}
\omega(\tau;\lambda_i) = 0\quad\Rightarrow\quad a_i(\tau)=0.
\end{equation}

\noindent
\textbf{P.3} The coefficients $a_i(\tau)$ are uniformly bounded in $\tau$, i.e.,
\begin{equation}\label{eq:gmls-ub}
\exists C>0\quad \mbox{s.t.}
\sum_{i=1}^{N} |a_i(\tau)| \le C\quad\forall \tau\in V^* .
\end{equation}

In what follows, we refer to \textbf{P.1} as \emph{P-reproduction} property, $\Lambda$ as the \emph{sampling} set, $P$ - as the \emph{reproduction property} space, and $a_i(\tau)$ - the basis functions. Following \cite{WendlandBook} we define the basis functions as the optimal solution of the following Quadratic Program (QP):
\begin{equation}\label{eq:gmls-qp}
\mbox{minimize}\,\,\frac12\sum_{i=1}^{N}\frac{a_i(\tau)^2}{\omega(\tau;\lambda_i)}
\quad
\mbox{subject to}
\quad
\sum_{i=1}^{N}a_i(\tau)\lambda_i(p_k) = \tau(p_k)\ \ k=1,2\ldots,Q.
\end{equation}
The objective in \eqref{eq:gmls-qp} enforces the local support of the basis, while the constraint enforces $P$-reproduction. 
The $P$-unisolvency condition \eqref{eq:P-uni} is sufficient  for \eqref{eq:gmls-qp} to have a unique minimizer that satisfies \textbf{P.1}--\textbf{P.3}; see \cite[Theorem 4.9, p.44]{WendlandBook}. With the notation
\begin{equation}\label{eq:GMLS-mats}
\begin{array}{c}
\displaystyle
\bm{a}(\tau)=(a_i(\tau))_{i=1}^{N}\in \mathbb{R}^N;\quad
\bm{r}(\tau)=(\tau(p_k))_{k=1}^{Q} \in \mathbb{R}^Q;\quad
\\[2ex]
\displaystyle
W(\tau)=\mbox{diag}(\omega(\tau;\lambda_i))\in \mathbb{R}^{N\times N};\quad
P(\Lambda) = (\lambda_{i}(p_j)) \in \mathbb{R}^{Q\times N}
\end{array}
\end{equation}
the QP \eqref{eq:gmls-qp} assumes an equivalent algebraic form
\begin{equation}\label{eq:gmls-qpalg}
\mbox{minimize}\,\,\frac12 \bm{a}(\tau)^T W(\tau)\bm{a}(\tau)
\quad
\mbox{subject to}
\quad
P(\Lambda)\bm{a}(\tau) = \bm{r}(\tau).
\end{equation}
It is easy to see that the optimal solution of \eqref{eq:gmls-qpalg}, resp \eqref{eq:gmls-qp} is given by
\begin{equation}\label{eq:gmls-opt}
\bm{a}(\tau) = W(\tau) P^T(\Lambda)\left(P(\Lambda)W(\tau)P^T(\Lambda)\right)^{-1}\bm{r}(\tau)
\end{equation}
and so, the GMLS approximation \eqref{eq:gmls} assumes the form
$$
\widetilde{\tau}(u) 
= \bm{\ell}^T(u)\bm{a}(\tau) = \bm{\ell}^T(u)W(\tau) P^T(\Lambda)\left(P(\Lambda)W(\tau)P^T(\Lambda)\right)^{-1}\bm{r}(\tau),
$$
with $\bm{\ell}^T(u) = (\lambda_i(u))_{i=1}^N\in \mathbb{R}^N$.
Setting 
\begin{equation}\label{eq:bcoeff}
\bm{b}(\tau)^T = \bm{\ell}^T(u)W(\tau) P^T(\Lambda)\left(P(\Lambda)W(\tau)P^T(\Lambda)\right)^{-1}
\end{equation}
allows us to write the approximation of the target functional as
\begin{equation}\label{eq:alt-rep}
\widetilde{\tau}(u) 
= \bm{\ell}^T(u)\bm{a}(\tau) 
= \bm{b}(\tau)^T\bm{r}(\tau).
\end{equation}
In other words, we can switch from a representation of $\widetilde{\tau}(u)$ in terms of a linear combination of samples $\lambda_i(u)$ 
to a representation in terms of a linear combination of the values $\tau(p_k)$ of the target functional at the basis of $P$. 
It is straightforward to show that $\bm{b}(\tau)$ admits the following variational characterization
\begin{equation}\label{eq:gmls-b}
\bm{b}(\tau) 
= \argmin_{\bm{c}\in\mathbb{R}^Q}\frac12
\left(\bm{\ell}(u) - P^T(\Lambda)\bm{c}\right)^T W(\tau)\left(\bm{\ell}(u) - P^T(\Lambda)\bm{c}\right).
\end{equation}
The algebraic least-squares problem  \eqref{eq:gmls-b} is equivalent to finding the best weighted least-squares fit to the data $\{\lambda_i(u)\}$ out of the space $P$ in the sense that the function 
\begin{equation}\label{eq:gmls-pstar}
p^*:= \bm{b}(\tau)^T\bm{p} = \sum_{k=1}^{Q} b_k(\tau) p_k,
\end{equation}
where $\bm{p} = (p_1,\ldots,p_Q)$, satisfies 
\begin{equation}\label{eq:lsfit-b}
p^*= \argmin_{p\in P} 
\frac12\sum_{i=1}^N 
\left( \lambda_i(u) - \lambda_i(p)\right)^2 \omega(\tau;\lambda_i).
\end{equation}
This motivates the notation
\begin{equation}\label{eq:gmls-tau}
\widetilde{\tau}(p^*):= \bm{b}(\tau)^T\bm{r}(\tau) = \sum_{k=1}^{Q} b_k(\tau) \tau(p_k).
\end{equation}
for the GMLS approximation of the target functional $\tau$.

In some important settings, such as the approximation of functions and their derivatives from scattered point values, the characterization \eqref{eq:lsfit-b}, resp., \eqref{eq:gmls-b} and the alternative representation  $\widetilde{\tau}(u)=\widetilde{\tau}(p^*)=\bm{b}(\tau)^T\bm{r}(\tau)$ offer some computational advantages. We consider such settings and their specialization to our needs in the next section.

\section{GMLS approximation of vector fields and their divergence from directional components} \label{sec:GMLS-div}
In this section we apply the abstract GMLS framework to define a discrete divergence operator $DIV_i :F^i\mapsto C^i$. Our approach draws upon ideas of  Mirzaei et al \cite{mirzaei2011,Mirzaei_13_ARXIV} on GMLS approximations of scalar functions and their derivatives from point function values. Accordingly, we start with a brief summary of the relevant aspects of their work.

\subsection{GMLS approximation of scalar functions and their derivatives}\label{sec:GMLS-scalar}
Let $m>0$ be some fixed integer. We consider the GMLS framework in Section \ref{sec:GMLS} with 
$$
V = C^{m+1}(\Omega),\quad
P = P_m \subset V,\quad
\Lambda=\{\delta_{\bm{x}_1},\ldots,\delta_{\bm{x}_N}\},
\quad\mbox{and}\quad
\tau(u) = \delta_{\bm{x}}\circ D^\alpha,
$$
where $X=\{\bm{x}_i,\ldots,\bm{x}_N\}$ is a finite set of points in $\Omega$, $\bm{x}\in \Omega$ is a given point assumed to be distinct from the points in $X$ and $\alpha$ is a multi index. If $|\alpha|=0$ then $\tau(u) = \delta_{\bm{x}}$, i.e., the target functional reduces to point evaluation at $\bm{x}$.

We assume that $X$ is such that the associated sampling set $\Lambda$ is $P_m$-unisolvent. The following conditions are sufficient to ensure this property; see \cite[Theorem 4.7, p.41]{WendlandBook},\cite[Theorem 3.14, p.33]{WendlandBook}, \cite[Definition 4.4]{mirzaei2011}:
\smallskip

\noindent
\textbf{C.1} The region $\Omega \in \mathbb{R}^d$ is compact and satisfies an interior cone condition with radius $r>0$ and angle $\theta \in \left(0,\, \pi/2\right)$.
\smallskip

\noindent
\textbf{C.2}  The point set $X= \{\bm{x}_i\}$ is quasi-uniform with fill distance (see \cite[Definition 4.6, p.41]{WendlandBook} $h_{X,\Omega} \le h_0$, where $h_0 = r/C(\theta)$ and $C(\theta) = 16(1+\sin\theta)^2 m^2/(3 \sin^2 \theta)$.
\medskip

Because both the target functional $\tau$ and the sampling set functionals $\lambda_i$ depend on the spatial location we can choose a weight function according to 
\begin{equation}\label{eq:wght}
\omega(\tau;\lambda_i) = \omega(\delta_{\bm{x}}\circ D^\alpha;\delta_{\bm{x}_i}) := \Phi(||\bm{x}-\bm{x}_i ||)
\end{equation}
where $\Phi(x)$ is a compactly supported, positive, radially symmetric function. In particular, we assume that the weight function satisfies the following condition:
\medskip

\noindent
\textbf{C.3} The support of the weight function $\Phi$ is contained in $\mathcal B(0, \delta)$ and $\Phi > 0$ in $\mathcal B(0, \delta/2)$, where $\delta = 2 C(\theta) h_{X,\Omega}$; see \cite[p.36]{WendlandBook}.
\medskip

Conditions \textbf{C.1}--\textbf{C.3} are sufficient to ensure that properties \textbf{P.1}-\textbf{P.3} hold in the present context. Thus, in what follows we restrict attention to GMLS settings satisfying these conditions. 

Recall the equivalent GMLS approximations of the target functional in \eqref{eq:alt-rep}. We now explain why in the present setting the representation $\widetilde{\tau}(u)=\bm{b}(\tau)^T\bm{r}(\tau)$ could be more computationally efficient than the equivalent representation  $\widetilde{\tau}(u)= \bm{\ell}^T(u)\bm{a}(\tau) $ in terms of the basis functions. 

Indeed, according to \eqref{eq:bcoeff} the coefficient vector $\bm{b}(\tau)$ depends on $\tau$ through the matrix $W(\tau)$. However, the choice of weight function \eqref{eq:wght} makes $W(\tau)$ dependent on the point location $\bm{x}$ \emph{but not on the derivative component} $D^\alpha$. It follows that  $\bm{b}$  will likewise depend on $\bm{x}$ but not on $D^\alpha$. In other words, 
\begin{equation}\label{eq:GMLS-derivative}
\widetilde{\tau}(u) = \widetilde{\delta_{\bm{x}}\circ D^\alpha}(u) 
= \bm{b}(\bm{x})^T\bm{r}(\delta_{\bm{x}}\circ D^\alpha).
\end{equation}
In contrast, \eqref{eq:gmls-opt} reveals that the basis functions $\bm{a}(\tau)$ depend both on $\bm{x}$ through $W(\bm{x})$ \emph{and the derivative component $D^\alpha$} through $\bm{r}(\tau) = \bm{r}(\delta_{\bm{x}}\circ D^\alpha)$, i.e., 
\begin{equation}\label{eq:GMLS-derivative-alt}
\widetilde{\tau}(u) = \widetilde{\delta_{\bm{x}}\circ D^\alpha}(u) 
= \bm{\ell}^T(u)\bm{a}(\delta_{\bm{x}}\circ D^\alpha). 
\end{equation}
Thus, if  $\delta_{\bm{x}}\circ D^\alpha$ is required for multiple values of $\alpha$, the coefficient vector $\bm{b}(\bm{x})$ must  be computed only once, whereas the basis functions $\bm{a}(\delta_{\bm{x}}\circ D^\alpha)$ must be recomputed for every new value of $\alpha$. Inversion of $P(\Lambda)W(\tau)P^T(\Lambda)$ comprises the bulk of the computational cost in both cases. Using \eqref{eq:GMLS-derivative} requires a single inversion of this matrix, whereas \eqref{eq:GMLS-derivative-alt} formally requires a separate inversion for every new target functional. 

\begin{remark}\label{rem:efficient}
In practice the actual computational advantage of \eqref{eq:GMLS-derivative} over \eqref{eq:GMLS-derivative-alt} depends on multiple factors such as the number of target functionals, the size of $X$, and the polynomial degree $m$ of the reconstruction space $P_m$. If  $m$ is small one could in principle store  $\left(P(\Lambda)W(\tau)P^T(\Lambda)\right)^{-1}$ in factored form at every point $\bm{x}_i$ and apply it to multiple right hand sides given by $\bm{r}(\delta_{\bm{x}}\circ D^\alpha)$. In this case the cost of \eqref{eq:GMLS-derivative-alt} is comparable to that of \eqref{eq:GMLS-derivative}. However, for larger polynomial degrees and/or large number of points the storage requirements may be prohibitively large, making \eqref{eq:GMLS-derivative} a better computational option. 
\end{remark}

In light of Remark \ref{rem:efficient} in what follows we use exclusively the representation \eqref{eq:GMLS-derivative}. Its computation decouples into a solution of a weighted least-squares problem 
\begin{equation}\label{eq:gmls-opt-x}
\bm{b}(\bm{x}) 
= \argmin_{\bm{c}\in\mathbb{R}^Q}\frac12
\left(\bm{\ell}(u) - P^T(\Lambda)\bm{c}\right)^T W(\bm{x})\left(\bm{\ell}(u) - P^T(\Lambda)\bm{c}\right)
\end{equation}
for the coefficients $\bm{b}(\bm{x})$, which is independent of the derivative part $D^\alpha$, followed by an application of the target functional to the basis functions of $P=P_m$. 
\begin{remark}\label{rem:diffuse}
Analogous to the relationship between \eqref{eq:gmls-b} and \eqref{eq:lsfit-b}, 
the algebraic least-squares
problem \eqref{eq:gmls-opt-x}  is equivalent to finding the best weighted least-squares polynomial fit to a set of function values
\begin{equation}\label{eq:lsfit-x}
p^*(\bm{x}) = \argmin_{p\in P_m} 
\frac12\sum_{i=1}^N 
\left( u(\bm{x}_i) - p(\bm{x}_i)\right)^2 \Phi(||\bm{x}-\bm{x}_i ||)
\end{equation}
in the sense that, given a basis $\{p_1,\ldots,p_Q\}$ of $P_m$, there holds
\begin{equation}\label{eq:opt-poly}
p^*(\bm{x}) = \bm{b}(\bm{x})^T \bm{r}(\delta_{\bm{x}}) = \sum_{i=1}^{Q} b_i(\bm{x}) p_i(\bm{x}).
\end{equation}
Likewise, for $\tau=\delta_{\bm{x}}\circ D^\alpha u$ the abstract formula \eqref{eq:gmls-tau} specializes to
\begin{equation}\label{eq:diffuse}
\widetilde{\delta_{\bm{x}}\circ D^\alpha}(u) = 
\widetilde{\delta_{\bm{x}}\circ D^\alpha}(p^*):=
\bm{b}(\bm{x})^T\bm{r}(\delta_{\bm{x}}\circ D^\alpha)=
\sum_{i=1}^{Q} b_i(\bm{x}) D^\alpha p_i(\bm{x}).
\end{equation}
The latter has the \emph{appearance} of an approximation to $D^\alpha p^*(\bm{x})$ computed by neglecting the dependence of the coefficients $b_i$ on the spatial location. This is why  \eqref{eq:diffuse} has been often referred to as a ``diffuse derivative'' \cite{Nayroles_92_CM} and deemed inferior to the true derivative $D^\alpha p^*(\bm{x})$. This misconception had been first pointed out by \cite{mirzaei2011}, which also provided a correct interpretation of \eqref{eq:diffuse} using the GMLS framework. Following their suggestion we refer to $\widetilde{\delta_{\bm{x}}\circ D^\alpha}(u)$ as the ``GMLS derivative approximation'', or simply the ``GMLS derivative'', instead of ``diffuse derivative''.
\end{remark}

To simplify notation we adopt the more compact symbol
$$
\widetilde{D^\alpha}u(\bm{x}) := \widetilde{\delta_{\bm{x}}\circ D^\alpha}(u)
$$
for the GMLS derivative. Assuming that conditions \textbf{C.1}--\textbf{C.3} hold Mirzaei et al \cite[Corollary 4.13]{mirzaei2011} prove that this derivative satisfies the following error bound:
\begin{equation}\label{eq:GMLS-error}
| \widetilde{D^\alpha}u(\bm{x}) - D^\alpha u(\bm{x})| \leq C h_{X,\Omega}^{m+1-|\alpha|} \|u\|_{C^{m+1}}.
\end{equation}

\subsection{A GMLS divergence operator} \label{sec:GMLS-vector}
%
We now focus on the application of the GMLS approach in Section \ref{sec:GMLS-scalar}  to define a discrete divergence operator $DIV_i :F^i\mapsto C^i$. This operator should map the normal vector field components defined on the local virtual cell faces $\bm{f}_{ij}$ to a single value on the local virtual cell $C_i$. Thus, to discuss the construction of $DIV_i$ it suffices to consider an arbitrary point $\bm{x}_i\in\Omega_h$ and its associated local primal-dual grid on $N^\varepsilon_i$. 

Owing to the duality relationship between $F_i$ and $E_i$, the normal components of $\bm{u}$ across the virtual faces $\bm{f}_{ij}$ coincide with the tangent components $u_{ij}=\bm{u}_{ij}\cdot\bm{t}_{ij}$ of this field along
the local edges $\bm{e}_{ij}$. 
Thus,  the task of defining $DIV_i$ essentially boils down to the task of approximating $\nabla\cdot\bm{u}(\bm{x}_i)$ from this data. We also note that using a set of non-normalized tangent vectors to the edges results in an equivalent set of degrees of freedom. In particular, the choice
\begin{equation}\label{eq:alt-tan}
u_{ij} = \bm{u}(\bm{x}_{ij})\cdot 2\bm{m}_{ij} = \bm{u}(\bm{x}_{ij})\cdot 2(\bm{x_{ij}}-\bm{x}_i)
\end{equation}
brings about some simplifications in the subsequent formulas and so, from now on we assume that the elements of $E^i$ are defined by \eqref{eq:alt-tan}.

Given a vector field $\bm{u}\in (C^{m+1}(\Omega))^d$ we call the scalar function
\begin{equation}\label{eq:radial}
u_{i\rightarrow}(\bm{x}) =  \bm{u}(\bm{x})\cdot 2(\bm{x}-\bm{x}_i)
\end{equation}
the radial component of $\bm{u}$ in the direction of $\bm{x}_i$. Clearly, $u_{i\rightarrow}\in C^{m+1}(\Omega)$,
\begin{equation}\label{eq:prop}
u_{i\rightarrow}(\bm{x}_{ij}) = \bm{u}(\bm{x}_{ij})\cdot 2\bm{m_{ij}} = u_{ij}
\quad\mbox{and}\quad u_{i\rightarrow}(\bm{x}_i) = 0.
\end{equation}

The following result shows that the vector field $\bm{u}$ and its divergence at $\bm{x}_i$ can be both expressed in terms of the derivatives of its radial component function.
\begin{lemma}\label{lem:radial-lap}
Let $\bm{u}\in (C^{m+1}(\Omega))^d$ and $u_{i\rightarrow}\in C^{m+1}(\Omega)$ be its radial component function \eqref{eq:radial}. There holds
\begin{equation}\label{eq:radial-lap}
\bm{u}(\bm{x}_i) = \frac12\nabla u_{i\rightarrow}(\bm{x}_i)
\quad\mbox{and}\quad
\nabla\cdot\bm{u}(\bm{x}_i) =\frac14 \Delta u_{i\rightarrow}(\bm{x}_i) .
\end{equation}
\end{lemma}
\emph{Proof.} 

Taking the gradient of the radial component function gives
$$
\nabla u_{i\rightarrow}(\bm{x}) = 2 \nabla \bm{u}(\bm{x})^T (\bm{x}-\bm{x}_i) + 2 \bm{u}(\bm{x}).
$$
Taking the divergence of this identity gives
$$
\begin{array}{rcl}
\nabla \cdot \nabla u_{i\rightarrow}(\bm{x}) &=& 2 \big(\nabla \cdot \nabla \bm{u} (\bm{x}) \big)\cdot (\bm{x}-\bm{x}_i) + 2 \nabla \bm{u}(\bm{x}) :\nabla \bm{x} + 2 \nabla \cdot \bm{u}(\bm{x}) \\[0.5ex]
& =&  2 \big(\nabla \cdot \nabla \bm{u} (\bm{x})\big) \cdot (\bm{x}-\bm{x}_i) + 4 \nabla \cdot \bm{u}(\bm{x})
\end{array}
$$

Setting $\bm{x}=\bm{x}_i$ completes the proof.
\endproof
\medskip

\begin{remark}

The line integral of $\bm{u}$ along the line $\bm{l}(t) = \bm{x_i} + 2t(\bm{x}- \bm{x_i})$, $t \in [0,1]$ provides an alternative definition of the radial component function
\begin{equation}\label{eq:integral}
u_{i\rightarrow}(\bm{x})  = \int_{\bm{x_i}}^{2\bm{x}- \bm{x_i}} \bm{u} \cdot d\bm{l} = \int_0^1 \bm{u}\left(\bm{l}(t)\right) \cdot 2(\bm{x}- \bm{x_i}) dt .
\end{equation}
Inserting the Taylor expansion  
$$
\bm{u}(\bm{l}(t)) = \bm{u}(\bm{x}_i) + 2 t \nabla \bm{u}(\bm{x}_i) (\bm{x}-\bm{x}_i) + o(\|\bm{x}-\bm{x}_i\|)
$$
of $\bm{u}(\bm{l}(t))$ about $t=0$ into  \eqref{eq:integral} and integrating the result gives
$$
u_{i\rightarrow}(\bm{x}) = \big( \bm{u}(\bm{x}_i) + \nabla \bm{u}(\bm{x}_i) (\bm{x}-\bm{x}_i) )\big) \cdot  2(\bm{x}- \bm{x_i}) + o(\|\bm{x}-\bm{x}_i\|^2,
$$
while the  expansion $\bm{u}(\bm{x}) = \bm{u}(\bm{x}_i) + \nabla \bm{u}(\bm{x}_i) (\bm{x}-\bm{x}_i) + o(\|\bm{x}-\bm{x}_i\|^2)$ implies that
$$
u_{i\rightarrow}(\bm{x}) = \bm{u}(\bm{x})\cdot 2(\bm{x}-\bm{x}_i) + o(\|\bm{x}-\bm{x}_i\|^2)
$$
Therefore, Lemma \ref{lem:radial-lap} continues to hold with the alternative definition \eqref{eq:integral}.

\end{remark}
\smallskip

Lemma \ref{lem:radial-lap} provides the foundation for a GMLS approximation of a vector field and its divergence from scattered directional components.  Specifically, the lemma asserts that we can define a GMLS approximation of $\bm{u}(\bm{x}_i)$ and $\nabla\cdot\bm{u}({\bm{x}_i})$ as 
\begin{equation}\label{eq:div-mls}
\widetilde{\bm{u}}(\bm{x}_i)=
\widetilde{\delta}_{\bm{x}_i}(\bm{u}):= \frac12 \widetilde{\nabla}u_{i\rightarrow}(\bm{x}_i)
\quad\mbox{and}\quad
\widetilde{\nabla\cdot}\bm{u}(\bm{x}_i):=\frac14 \widetilde{\Delta}u(\bm{x}_i) ,
\end{equation}
respectively, where $\widetilde{\nabla}u_{i\rightarrow}(\bm{x}_i)$ and $\widetilde{\Delta}u_{i\rightarrow}(\bm{x}_i)$ are the GMLS gradient and Laplacian of the radial component function, respectively. 

Our main idea now is to specialize the GMLS approach in Section \ref{sec:GMLS-scalar} to the radial component function in such a way that the GMLS divergence  becomes a mapping $F^i\mapsto C^i$. 
Such a specialization enables us to use  $\widetilde{\nabla\cdot}\bm{u}(\bm{x}_i)$ to define the discrete divergence operator  $DIV_i$ in our staggered scheme.

In light of \eqref{eq:prop}, i.e., $u_{i\rightarrow}(\bm{x}_i)=0$, the appropriate setting for this specialization is given by the spaces
$$
V=C^{m+1}_i := \{ v\in C^{m+1}(\Omega)\,|\, v(\bm{x}_i) = 0\}
\quad\mbox{and}\quad
P= P^i_m:=\{ p\in P_m \,|\, p(\bm{x}_i) = 0\}
$$
while
\begin{equation}\label{eq:sample-div}
\Lambda_i=\{\delta_{\bm{x}_{ij}}\,|\, \bm{x}_{ij}\in M_i\}
\quad\mbox{and}\quad
\tau_i(u) = \delta_{\bm{x}_i}\circ D^\alpha
\end{equation}
define the set of the sampling functionals and the target functionals, respectively. As in section \ref{sec:GMLS-scalar}, this allows us to use the weight function \eqref{eq:wght}.

We proceed to establish a local $P^i_m$-reproduction property. 
\begin{lemma}
\label{polynomialRecAlt}
Assume that \textbf{C.1}--\textbf{C.3} hold for $X=M_i$, i.e., the set of midpoints in the neighborhood $N^\varepsilon_i$. Let $\Lambda_i$ and $\tau_i$ be the functionals in  \eqref{eq:sample-div}. Then, there exist coefficients ${a}_{ij,\alpha}(\bm{x}_i)={a}_{ij}(\delta_{\bm{x}_i}\circ D^\alpha)$ such that properties \textbf{P.1}--\textbf{P.3} hold with $P=P^i_m$. In particular, there exist a constant $C_{1,\alpha}>0$
such that ${a}_{ij,\alpha}(\bm{x}_i)$ satisfy
\begin{equation} \label{bds}
\sum_{j=1}^N |{a}_{ij,\alpha}(\bm{x}_i)| \le C_{1,\alpha} 
\end{equation}
and ${a}_{ij,\alpha}(\bm{x}_i) = 0$, if $\|\bm{x}_i-\bm{x}_{ij}\| > 2 C(\theta) h_{X, \Omega}$, where $h_{X, \Omega}$ and $C(\theta)$ are as in Condition \textbf{C.2}.
The equivalent problem \eqref{eq:lsfit-b} has a unique solution $\bm{b}(\bm{x})$, which defines a GMLS derivative according to \eqref{eq:diffuse}.
\end{lemma}
\smallskip

{\em Proof}.
Conditions \textbf{C.1}--\textbf{C.3} ensure that the midpoints $M_i$ are $P_m$-unisolvent and so, the sampling set $\Lambda\mauro{_i}$ is also $P_m$-unisolvent. Since $P^i_m\subset P_m$ it follows that $\Lambda\mauro{_i}$ is also unisolvent for $P^i_m$. The rest of the proof follows the arguments in \cite{WendlandBook} and \cite{mirzaei2011}.
\endproof
\medskip

The next Lemma estimates the errors of the specialized GMLS approximation.
\smallskip

\begin{lemma}
\label{polynomialRec2}
Assume that the hypotheses of Lemma \ref{polynomialRecAlt} and let $\Omega^* = \overline{\mathcal B({\bm{x}_i},\delta)}$. For any $u \in C^{m+1}(\Omega^*)\cap C^{m+1}_i(\Omega)$  there holds
$$
\left \|  \widetilde{D}^\alpha u ({\bm{x}_i}) - D^{\alpha} u({\bm{x}_i}) \right \| \le C h_{X,\Omega}^{m+1-|\alpha|} \| u\|_{C^{m+1}(\Omega^*)}.
$$
\end{lemma}
{\em Proof}. We adapt the proof of Theorem 4.3 in \cite{mirzaei2011}. Adding and subtracting an arbitrary $p\in P^i_m$ and using the $P^i_m$ reproduction property yields the bound
\begin{equation} \label{bound1}
\begin{array}{l}
\left \|  D^{\alpha} u({\bm{x}_i}) - \widetilde{D}^\alpha u({\bm{x}_i}) \right \| \le \left \|  D^{\alpha} u ({\bm{x}_i}) - D^{\alpha} p({\bm{x}_i}) \right \|+ \left \|   D^{\alpha} p({\bm{x}_i}) - \widetilde{D}^\alpha u({\bm{x}_i}) \right \|\\[2ex]
\quad = \displaystyle \left\|  D^{\alpha} u({\bm{x}_i}) - D^{\alpha} p({\bm{x}_i}) \right \| + 
\Big \|  \sum_{ij} a_{ij,\alpha}({\bm{x}_i}) p(\bm{x}_{ij})  - \sum_{ij} a_{ij,\alpha}({\bm{x}_i}) u({\bm{x}_{ij}})  \Big \| \\[2ex]
\quad \le \left \|  D^{\alpha} u ({\bm{x}_i}) - D^{\alpha} p({\bm{x}_i}) \right \| + \|p - u \|_{L^{\infty}(\Omega^*)} \sum_{ij} | a_{ij,\alpha}({\bm{x}_i})|
\end{array}
\end{equation}
Because $u({\bm{x}_i})=0$, its Taylor polynomial $p^i_m$ of order $m$ about ${\bm{x}_i}$ belongs to $P^i_m$. Setting $p=p^i_m$ in \eqref{bound1} and using that 
\begin{equation} \label{bound2}
\left \|  D^{\alpha} u - D^{\alpha} p^i_m \right \|_{L^{\infty}(\Omega^*)}  \le C h_{X,\Omega}^{m+1-|\alpha|} \| u \|_{C^{m+1}(\Omega^*)}.
\end{equation}
completes the proof. 
\endproof
\bigskip

\begin{corollary}\label{cor:error}
Assume the hypotheses and notation of Lemmas \ref{polynomialRecAlt}--\ref{polynomialRec2}.
Let $\widetilde{\bm{u}}(\bm{x}_i)$ and $\widetilde{\nabla\cdot}\bm{u}(\bm{x}_i)$ be the GMLS approximations of 
$\bm{u}\in (C^{m+1}(\Omega))^d$ and its divergence defined in \eqref{eq:div-mls}. Then,
$$
\big \|  \bm{u}(\bm{x}_i) - \widetilde{\bm{u}}(\bm{x}_i) \big \| \le C h_{X,\Omega}^{m} \| \bm{u}\|_{C^{m+1}(\Omega^*)}.
$$
and
$$
\big \|  \nabla\cdot\bm{u}(\bm{x}_i) - \widetilde{\nabla\cdot}\bm{u}(\bm{x}_i) \big \| \le C h_{X,\Omega}^{m-1} \| \bm{u}\|_{C^{m+1}(\Omega^*)}.
$$
\end{corollary}

\begin{remark}
Although Lemma \ref{polynomialRecAlt} holds for a generic point $\bm{x}$, the error estimate in Lemma \ref{polynomialRec2} is valid only at $\bm{x}_i$.  The proof of the latter requires the Taylor polynomial of $u$ about the point where the error is being estimated to belong in $P^i_m$, i.e., to vanish at $\bm{x}_i$. In general, if $p_m$ is the Taylor polynomial of $u$ about some other point $\bar{\bm{x}}$ there is no guarantee that $p_m\in P^i_m$. Nonetheless, this does not present a problem, since we are only interested in estimates of the derivatives of $u$ at $\bm{x}_i$.
\end{remark}

In summary, given a vector field $\bm{u}\in (C^{m+1}(\Omega))^d$ and a point $\bm{x}_i$ with neighborhood $N^\varepsilon_i$, our GMLS approach for the approximation of $\bm{u}$ and its divergence at $\bm{x}_i$ comprises the following two steps:
\medskip

\begin{shaded}
\noindent
1.  Solve the weighted least-squares problem
\begin{equation} \label{auxiliaryproblem}
p^*(\bm{x}_i) = \argmin\limits_{p\in P_m^i}  \Big\{ \sum_{j=1}^N \left[ \bm{u}(\bm{x}_{ij})\cdot 2\bm{m}_{ij} - \, p(\bm{x}_{ij}) \right]^2 \Phi(\|\bm{x}_i-\bm{x}_{ij}\|) \Big\}.
\end{equation}

\noindent
2.  Let $\bm{b}(\bm{x}_i)$ denote the coefficient vector of $p^*$ relative to a basis $\{p_k^i\}$ of $P_m^i$. Set
\begin{equation}\label{eq:GMLS-approx}
\begin{array}{rl}
\widetilde{\bm{u}}(\bm{x}_i)  =&\!\!\!
\displaystyle
\frac12\widetilde{\nabla}p^*({\bm x}_i) :=\frac12\sum_{k=1}^Q b_k(\bm{x}_i)\nabla p_k^i(\bm{x}_i), \\[3ex]
\widetilde{\nabla\cdot}\bm{u}(\bm{x}_i)= &\!\!\!
\displaystyle\frac14\widetilde{\Delta}p^*({\bm x}_i):=
\frac14\sum_{k=1}^Q b_k(\bm{x}_i)\Delta p_k^i(\bm{x}_i).
\end{array}
\end{equation}
\end{shaded}

\section{Formulation of the staggered GMLS discretization method}\label{sec:define}

The duality of the local geometric structures and the ensuing isomorphism of $E^i$ and $V^i$ allow us to discretize the second-order operator in \eqref{eq:model} by ``chaining''  the discrete divergence and gradient operators at every point of $\Omega_h$.  This observation is at the heart of our staggered GMLS method. 
For clarity we present the method for $\mu=I$ and Dirichlet boundary conditions, i.e., $\Gamma_D = \Gamma$. 
Section \ref{sec:implement} briefly discusses implementation of the method for a general $\mu$ and Neumann boundary conditions.

Let $\Omega_h$ be a meshless discretization of the computational domain $\Omega$ such that the points in $\Omega_h$ satisfy conditions \textbf{C.1}--\textbf{C.2}. The global meshless approximation space   
$$
V^h=\bigcup_{i=1}^{N} V^i 
$$
is the union of the local vertex spaces associated with the points in $\Omega_h$. The gradient of a discrete function $\phi^h\in V^h$ is simply the application of $GRAD_i$ to the restriction of $\phi^h$ to $V^i$. 

As in many other discretization methods we choose to impose the Dirichlet boundary condition on the approximating space by setting $\phi(\bm{x}_i) = u(\bm{x}_i)$ for all boundary particles. Thus, our new staggered GMLS discretization of  \eqref{eq:model} reads: seek $\phi^h\in V^h$ such that
\begin{equation}\label{eq:order2}
\begin{array}{rcl}
DIV_i \circ GRAD_i (\phi^h) & = & f(\bm{x}_i)
\quad\forall i\in \breve{\Omega}_h \\[1ex]
\phi_i & = & u(\bm{x}_i) \quad\forall i\in \Gamma_h .
\end{array}
\end{equation}
To assemble the discrete problem  \eqref{eq:order2} at every interior point ${\bm x}_i$  we first compute the topological gradient of the local discrete field and then apply the GMLS divergence operator. The second step follows the procedure outlined in the gray box in Section \ref{sec:GMLS-vector}, i.e., we solve the weighted least-squares problem 
\begin{equation} \label{eq:stagMLS}
p^*(\bm{x}_i) = \argmin\limits_{p\in P_m^i}  \Big\{ \sum_{j=1}^N \left[ GRAD_i (\phi^h) - \, p(\bm{x}_{ij}) \right]^2 \Phi(\|\bm{x}_i-\bm{x}_{ij}\|) \Big\}
\end{equation}
where $\Phi(\cdot)$ satisfies condition \textbf{C.3} and then set
$$
DIV_i\circ GRAD_i (\phi^h)  = \frac14 \widetilde{\Delta}p^*({\bm x}_i).
$$

In the next section we show that the GMLS gradient of $p^*$ is an accurate approximation of the flux at the virtual cell center, i.e.,
$$
\bm{u}(\bm{x}_i)=-\nabla \phi (\bm{x}_i) \approx \frac12 \widetilde{\nabla}p^*(\bm{x}_i).
$$
This result is useful in settings that also require the flux values.  
\medskip

\begin{remark}\label{rem:compatible}
The scheme \eqref{eq:order2} is equivalent\footnote{The two forms of the staggered GMLS are equivalent in the sense that both produce the exact same discrete approximation $\phi^h\in V^h$.} to the following ``mixed'' system of equations: seek $\bm{u}^h\in F^h$ and $\phi^h\in V^h$ such that $\phi_i=u(\bm{x}_i)$ for all $i\in\Gamma_h$ and 
\begin{equation}\label{eq:mixed}
\left\{
\begin{array}{rcl}
DIV_i(\bm{u}^h) &=& f(\bm{x}_i) \\[1ex]
\bm{u}^h +  GRAD_i (\phi^h) &=&0
\end{array}
\right.
\quad\forall i\in \breve{\Omega}_h ,
\end{equation}
where $F^h=\cup_{i=1}^N F^i$. Note that all metric attributes of the staggered GMLS scheme are subsumed in the first equation, whereas the second equation is purely topological, i.e., it depends only on the local connectivity graph of $N^\varepsilon_i$ and not on the physical locations of its points. This structure mirrors the structure of many compatible discretization methods for \eqref{eq:model}, which also combine a metric-dependent equation with a topological relation \cite{Bochev_06a_IMA}. 
A more rigorous qualification of the compatibility properties of the new scheme is beyond the scope of this work and will be pursued in a forthcoming article. However, it is worth pointing out that numerical results in Section  \ref{sec:results} confirm that solutions of \eqref{eq:order2} for problems with discontinuous material properties are indeed qualitatively very similar to solutions of compatible discretization methods for \eqref{eq:model}. 
\end{remark}


\section{Accuracy of the staggered GMLS method}
\label{sec:props}
In this section we study the polynomial consistency of the staggered GMLS method \eqref{eq:order2}.

\begin{theorem}
\label{gradTheorem}
Assume that conditions \textbf{C.1}--\textbf{C.3} hold. Let $\phi^h\in V^h$ be the interpolant of a function $\phi \in C^{m+1}(\Omega)$ and let $p^*({\bm x}_i)$ be the corresponding solution of \eqref{eq:stagMLS}. 
Then, the following error bounds hold:
\begin{align}
&\|\Delta \phi (\bm{x}_i)  - DIV_i \circ GRAD_i (\phi^h)  \|\le C h_{X,\Omega}^{m-1} \| \phi\|_{C^{m+1}(\Omega)} \\[1ex]
&\|\nabla \phi (\bm{x}_i) -\frac12 \widetilde{\nabla}{p}_i^*(\bm{x}_i) \|\le C h_{X,\Omega}^{m} \| \phi\|_{C^{m+1}(\Omega)}.
\end{align}
\end{theorem}

\emph{Proof.} Consider the auxiliary vector field
$$
\bm{u}(\bm{x}) = g(\bm{x})\frac{\bm{x}-\bm{x}_i}{2\|\bm{x}-\bm{x}_i\|^2},\quad\mbox{where}\quad
g(\bm{x})=\phi(2\bm{x}-\bm{x}_i)-\phi(\bm{x}_i).
$$
The radial component of $\bm{u}(\bm{x})$ equals\footnote{The alternative definition \eqref{eq:integral} with $\bm{u} = \nabla \phi$ yields the same result, i.e. $u_{i\rightarrow}(\bm{x})= g(\bm{x})$.} $g(\bm{x})$, i.e.,
$$
u_{i\rightarrow}(\bm{x})= \, \bm{u}(\bm{x})\cdot 2 (\bm{x}-\bm{x}_i) = g(\bm{x}),
$$
Furthermore, $g(\bm{x}_i)=0$, that is $g(\bm{x})\in C^{m+1}_i$, and 
$$
g(\bm{x}_{ij})=\phi(\bm{x}_j)-\phi(\bm{x}_i)=GRAD_i(\phi^h).
$$
It follows that the weighted least-squares problem \eqref{eq:stagMLS} in the formulation of our method is equivalent to the GMLS problem \eqref{auxiliaryproblem} for the auxiliary velocity field $\bm{u}(\bm{x})$.  Accordingly, the identities in \eqref{eq:GMLS-approx} specialize to
$$
\widetilde{\bm{u}}(\bm{x}_i)  =
\displaystyle\frac12\widetilde{\nabla}p^*({\bm x}_i)  
\quad\mbox{and}\quad
\widetilde{\nabla\cdot}\bm{u}(\bm{x}_i)= 
\frac14\widetilde{\Delta}p^*({\bm x}_i) =DIV_i\circ GRAD_i (\phi^h) ,
$$
respectively. 
Using the identities
$$
\nabla g(\bm{x}) = 2 \nabla \phi(\bm{x})\quad\mbox{and}\quad 
\Delta g(\bm{x}) = 4 \Delta \phi(\bm{x}) \,,
$$
and the fact that $g(\bm{x})$ is the radial component of $\bm{u}(\bm{x})$ implies that
$$
\bm{u}(\bm{x}_i) = \frac12 \nabla g(\bm{x}_i) = \nabla \phi(\bm{x}_i) 
\quad\mbox{and}\quad
\nabla\cdot\bm{u}(\bm{x}_i) = \frac14 \Delta g(\bm{x}_i)= \Delta \phi(\bm{x}_i),
$$
respectively. Therefore,
$$
\nabla \phi (\bm{x}_i) -\frac12 \widetilde{\nabla}{p}_i^*(\bm{x}_i)
=
\bm{u}(\bm{x}_i)  -\widetilde{\bm{u}}(\bm{x}_i),
$$
and
$$
\Delta \phi (\bm{x}_i)  - DIV_i\circ GRAD_i (\phi^h) 
=
\nabla\cdot\bm{u}(\bm{x}_i) - \widetilde{\nabla\cdot}\bm{u}(\bm{x}_i),
$$
respectively. Application of the error bounds in Corollary \ref{cor:error} grants the proof.
\endproof

\section{Extensions and implementation details}\label{sec:implement}
This section briefly reviews the  treatment of a variable material property tensor $\mu$ as well as the imposition of the Neumann boundary condition  in the staggered GMLS. We also comment on some practical aspects of the implementation used in our numerical studies, including choice of polynomial basis and solution of the local optimization problems.

Extension of \eqref{eq:order2} to a variable $\mu$ is straightforward. Assuming that $\mu$ is available at all points in $\Omega_h$ we solve a modified version of \eqref{eq:stagMLS} given by 
\begin{equation} \label{eq:stagMLSmu}
p^*(\bm{x}_i) = \argmin\limits_{p\in P_m^i}  \Big\{ \sum_{j=1}^N \left[ \mu_{ij}GRAD_i (\phi^h) - \, p(\bm{x}_{ij}) \right]^2 \Phi(\|\bm{x}_i-\bm{x}_{ij}\|) \Big\}
\end{equation}
where $GRAD_i$ is the topological gradient defined in \eqref{eq:grad} and $\mu_{ij} = (\mu(\bm{x}_i)+\mu(\bm{x}_j))/2$. Then we approximate the second-order operator in \eqref{eq:model} according to
$$
\nabla\cdot\mu\nabla \phi(\bm{x}_i) \approx DIV_i\circ \mu GRAD_i (\phi^h) := \frac14 \widetilde{\Delta} p^*(\bm{x}_i).
$$

Implementation of the Neumann condition $\partial_n \phi = g$ in the staggered GMLS scheme is somewhat more involved. 
Recall that in the case of Dirichlet conditions we approximate the second-order operator $\nabla\cdot\mu\nabla \phi(\bm{x}_i)$ only at the interior points and impose the boundary condition simply by setting  $\phi_i=u(\bm{x}_i)$ at all $\bm{x}_i\in \Gamma_D$. 
To impose the Neumann condition, in addition to interior points, we also approximate $\nabla\cdot\mu\nabla \phi(\bm{x}_i)$ at all Neumann points $\bm{x}_i\in \Gamma_N$ and incorporate the boundary condition in the definition of the discrete second-order operator.  The latter is accomplished by constraining  \eqref{eq:stagMLSmu} with a GMLS approximation of the equation $\partial_n \phi = g$. Taking into account that $\nabla\phi(\bm{x}_i)\approx 1/2 \widetilde{\nabla}{p}^*(\bm{x}_i)$, this approximation is given by
$$
\widetilde{\partial_n}p(\bm{x}_i) = 2 g(\bm{x}_i)
$$
and results in the following constrained version of \eqref{eq:stagMLSmu}:
\begin{equation} \label{eq:stagMLSmu-constr}
p^*_{\partial n}(\bm{x}_i) = \argmin\limits_{p\in P_m^i; \widetilde{\partial_n}p(\bm{x}_i) = 2 g(\bm{x}_i)}  \Big\{ \sum_{j=1}^N \left[ \mu_{ij}GRAD_i (\phi^h) - \, p(\bm{x}_{ij}) \right]^2 \Phi(\|\bm{x}_i-\bm{x}_{ij}\|) \Big\}.
\end{equation}
Solving \eqref{eq:stagMLSmu-constr} instead of \eqref{eq:stagMLSmu} effectively replaces the operator $DIV_i$ with a modified version $DIV_{\partial n,i}$ which naturally incorporates the Neumann boundary condition. We then impose this condition by approximating the second-order operator at all points on the Neumann boundary according to 
$$
\nabla\cdot\mu\nabla \phi(\bm{x}_i) \approx DIV_{\partial n,i}\circ \mu GRAD_i (\phi^h) 
:= \frac14 \widetilde{\Delta} p^*_{\partial n}(\bm{x}_i).
$$ 
where $p^*_{\partial n}$ solves \eqref{eq:stagMLSmu-constr}.

\subsection{Practical aspects of the implementation}
In the general case of a variable $\mu$ and mixed boundary conditions the staggered GMLS scheme comprises the following system of algebraic equations:
\begin{equation}\label{eq:GMLS-all}
\begin{array}{rcl}
 DIV_i\circ \mu GRAD_i (\phi^h) &=& f_i\quad \forall \bm{x}_i\in  \breve{\Omega}_h \\[1ex]
 DIV_{\partial n,i}\circ \mu GRAD_i (\phi^h) &=&  f_i\quad \forall \bm{x}_i\in \Gamma_N \\[1ex]
 \phi_i &=& u_i \quad \forall \bm{x}_i\in \Gamma_D
\end{array}.
\end{equation}
Assembly of this system requires solution of \eqref{eq:stagMLSmu} for the discretization at all interior points and solution of the constrained problem \eqref{eq:stagMLSmu-constr} to compute the second-order discrete operator at all Neumann points.  The actual form of these problems depends on the choice of a basis for the polynomial reconstruction. In our implementation we choose the scaled and shifted Taylor monomials, i.e.,  
\begin{equation}\label{eq:taylor}
P_m^i = span \left\{ p_\alpha^i \right \}_{|\alpha | \le m}, \quad p_\alpha^i = 
\frac{1}{\alpha!} \left(\frac{\bm{x}-\bm{x}_i}{\epsilon}\right)^\alpha.
\end{equation}
In the present context the matrices $P$, $W$ and the vector $\bm{\ell}$ in Section \ref{sec:GMLS} specialize to 
$$
W_i=\mbox{diag}( \Phi (||\bm{x}_{ij} - \bm{x}_i||))_{\bm{x}_{ij}\in M_i}, \ 
P_i=\left(p_\alpha^i(\bm{x}_{ij})\right)_{\bm{x}_{ij}\in M_i}^{|\alpha|\le m}, \
\bm{\ell}_i = \left(\mu_{ij} (\phi_j - \phi_i)\right)_{\bm{x}_{ij}\in M_i}.
$$
The optimal coefficient vector $\bm{b}_i$ for the solution of \eqref{eq:stagMLSmu} is given by \eqref{eq:gmls-b}, which specializes to
\begin{equation}\label{eq:opt-coeff}
\bm{b}^T_i = \bm{\ell}_i^T\,W\, P_i^T \left( P_iW_iP_i^T \right)^{-1}\,.
\end{equation}
As a result, the approximation of the second-order operator at all interior points assumes the form
$$
DIV_i\circ \mu GRAD_i (\phi^h) =\frac14 \bm{b}^T_i \bm{r}(\delta_{\bm{x}_i}\circ\Delta)
$$
where $\bm{r}(\delta_{\bm{x}_i}\circ\Delta) := \left(\Delta p_\alpha^i(\bm{x}_i)\right)$.


To solve the constrained problem  \eqref{eq:stagMLSmu-constr} we introduce a Lagrange multiplier $\lambda$ and seek the stationary point $\{p^*_{\partial n},\lambda^*\}$ of the Lagrangian functional 
$$
L(p,\lambda)= \sum_{j=1}^N \left[ \mu_{ij}GRAD_i (\phi^h) - \, p(\bm{x}_{ij}) \right]^2 \Phi(\|\bm{x}_i-\bm{x}_{ij}\|)
-\lambda\left(\widetilde{\partial_n}p(\bm{x}_i) - 2 g(\bm{x}_i\right).
$$ 
Let $\bm{c}^*_i = \left( \bm{b}_{\partial n,i},\lambda\right)$ be the coefficient vector of $\{p^*_{\partial n},\lambda^*\}$ and $\bm{r}(\delta_{\bm{x}_i}\circ \partial_n)= \left( \partial_n p_\alpha^i(\bm{x}_i)\right)$. It is not hard to see that 
the necessary optimality condition for the stationary point is equivalent to a linear algebraic system $K_i\bm{c}_i = \bm{d}_i$ for the optimal coefficient vector  $\bm{c}^*_i $, where
\begin{equation}
\label{lagrange1}
K_i =\left[
\begin{array}{cc}
P_iW_iP_i^T & \bm{r}(\delta_{\bm{x}_i}\circ \partial_n) \\[1ex]
\bm{r}(\delta_{\bm{x}_i}\circ \partial_n)^T & 0
\end{array}\right]
\quad\mbox{and}\quad
\bm{d}_i = \left[
\begin{array}{c}
P_iW_i\bm{\ell}_i\\ 2 g_i
\end{array}
\right].
\end{equation}
After solving \eqref{lagrange1} and extracting $ \bm{b}_{\partial n,i}$ we obtain the approximation 
$$
DIV_{\partial n,i}\circ \mu GRAD_i (\phi^h) = \frac14 \bm{b}^T_{\partial n,i} \bm{r}(\delta_{\bm{x}_i}\circ\Delta).
$$
of the differential operator at all Neumann points.


Computation of the optimal coefficient vectors for both \eqref{eq:stagMLSmu} and \eqref{eq:stagMLSmu-constr} requires inversion of the matrix $M_i=P_iW_iP_i^T $, which resembles  the poorly conditioned Hilbert matrix. Because in this work we used polynomial degrees $m\leq 6$, the relatively small size of the resulting local problems enabled us to stably invert $M_i$ by using a standard LU decomposition. Larger polynomial degrees may require more advanced solvers and/or a different choice of basis polynomials.

Due to the strong form discretization of the PDE and the lack of symmetry in particle arrangement, the global matrix corresponding to \eqref{eq:GMLS-all} is in general not symmetric. 
Furthermore, a special care must be taken to ensure solvability for problems with a non-trivial null-space as in the case of pure Neumann boundary conditions. Possible options for handling the nullspace include pinning a degree of freedom, adding a global Lagrange multiplier to remove the singularity, or projecting out the kernel if a Krylov subspace method is used to solve the resulting equations. For more information about this subject we refer to
\cite{traskCMAME,BochevLehoucq,Yeckel_99_IJNMF} and the references cited therein.

\begin{remark}
Consider an interior point $\bm{x}_i\in\breve{\Omega}$. By collecting the terms multiplying  the elements of $\bm{\ell}_i= \left(\mu_{ij} (\phi_j - \phi_i)\right)$ we can write the discrete operator as
$$
DIV_i\circ \mu GRAD_i (\phi^h) 
= \sum_{j\in N^\varepsilon_i\cap\breve{\Omega}} \beta_{ij} (\phi_j - \phi_i)
+ \sum_{j\in N^\varepsilon_i\cap\Gamma_D} \beta_{ij} (u_j - \phi_i)
$$
for a suitable set of coefficients $\beta_{ij}$. Likewise, at all Neumann points the discrete operator can be written as
$$
DIV_{\partial n, i}\circ \mu GRAD_i (\phi^h) 
= \sum_{j\in N^\varepsilon_i} \beta_{ij} (\phi_j - \phi_i) + \gamma_i g_i
$$
for some suitable coefficient $\gamma_i$. Consequently,  the system \eqref{eq:GMLS-all} has the following equivalent form: 
\begin{equation}\label{eq:GMLS-all-FD}
\begin{array}{rcll}
\displaystyle
 \sum_{j\in N^\varepsilon_i\cap\breve{\Omega}} \beta_{ij} (\phi_j - \phi_i) -
 \sum_{j\in N^\varepsilon_i\cap\Gamma_D} \beta_{ij}  \phi_i
 &=&
 \displaystyle f_i - \sum_{j\in N^\varepsilon_i\cap\Gamma_D} \beta_{ij} u_j & \forall \bm{x}_i\in  \breve{\Omega}_h
  \\[1ex]
  \displaystyle
 \sum_{j\in N^\varepsilon_i} \beta_{ij} (\phi_j - \phi_i)  &=&  f_i - \gamma_i g_i & \forall \bm{x}_i\in \Gamma_N 
 \\[1ex]
 \phi_i &=& u_i & \forall \bm{x}_i\in \Gamma_D
\end{array}.
\end{equation}
Thus, the staggered GMLS method can be interpreted as a finite-difference-like formula for the Laplacian defined on an irregular stencil comprising the points in $N^\varepsilon_i$.
\end{remark}

\section{Numerical examples}\label{sec:results}

\subsection{Case setup}
To discretize the domain for all problems in this work, we first discretize the boundary of the domain with a length scale $dx$ and then discretize the interior of the domain with a Cartesian lattice with spacing $dx$. We then perturb each interior particle by a uniformly distributed random variable of magnitude $\eta$, and delete all particles that lie outside of the domain. The variable $\eta$ is used to demonstrate insensitivity of the results to particle anisotropy, and for all results presented is set to $\eta = 0.1 dx$. We characterize levels of discretization by the number of particles in each direction of the lattice $N = \frac{1}{dx}$.

Algebraic multigrid is used to efficiently solve the discretized asymmetric system of equations. While an in-depth discussion of the performance of this preconditioner is left for a future work, we briefly comment that a combination of GMRes and auxilliary space preconditioning \cite{jinchao} provides an efficient means of solution with computational effort scaling nearly linearly with the degrees of freedom in the system.

\subsection{Elliptic problems with smooth coefficients}
\begin{figure}[h!]
  \centering
  \includegraphics[width=0.45\textwidth]{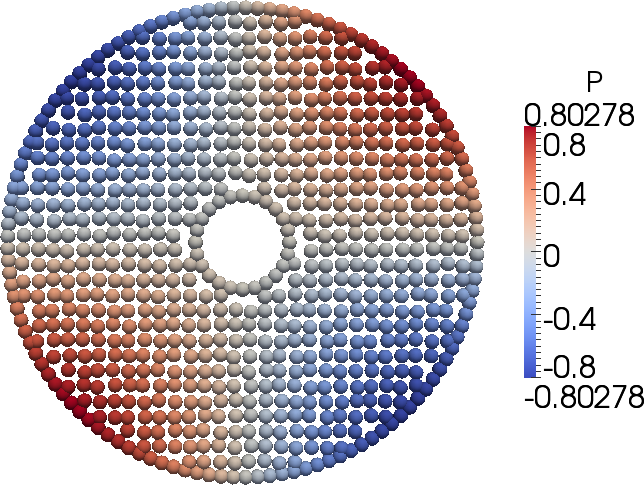}
  \includegraphics[width=0.4\textwidth]{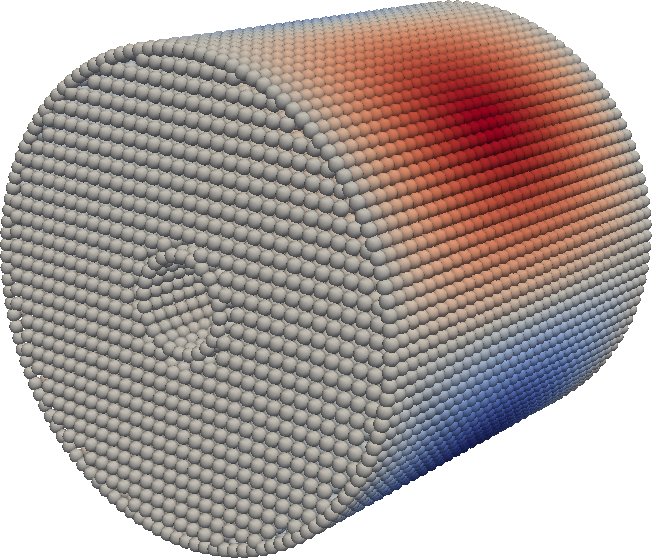}
\caption{Typical particle distribution with $N=32$ for 2D annular disk geometry (left) and 3D extrusion (right).}
  \label{manufacturedgeometry}
\end{figure}

We first demonstrate the convergence rate of the method by using smooth manufactured solutions and  considering both an annular 2D domain with outer radius $R_o = \frac{\pi}{2}$ and inner radius $R_i = \frac{\pi}{4}$, and a 3D domain formed by extruding the 2D geometry in the z-direction to form a cylinder with height $\pi$; see Figure \ref{manufacturedgeometry}. 
We set the exact solutions to
$$
\phi^{2D}_{ex} = \sin x \sin y
\quad\mbox{and}\quad
\phi^{3D}_{ex} = \sin x \sin y \sin z
$$
in two and three-dimensions, respectively. Substitution of the exact solutions in the model problem \eqref{eq:model} defines the source term and the boundary data. We then solve the discrete problem \eqref{eq:GMLS-all} with either Dirichlet or Neumann boundary conditions. 
As no mesh is available to quantify the error, we measure convergence in a root mean square sense by defining the following norm
\begin{equation}
\| \phi^h \|_{l_2} =\left(\frac1N\sum_{i\in\Omega_h}\phi_i^2\right)^{1/2}
\end{equation}
and semi-norm
\begin{equation}
|\phi^h |_{h^1} = \left(\frac1N\sum_{i\in\Omega_h}(GRAD_i\phi^h)^2\right)^{1/2}
\end{equation}

We first compare the convergence rates of the staggered GMLS scheme with those of a \emph{collocated} one obtained by using the standard GMLS derivative approximation of $\nabla\cdot\mu\nabla\phi$ \cite{mirzaei2011}. In both cases we use the same Taylor basis \eqref{eq:taylor} and kernel functions. 
\begin{figure}[h!]
  \centering
  \includegraphics[width=0.45\textwidth]{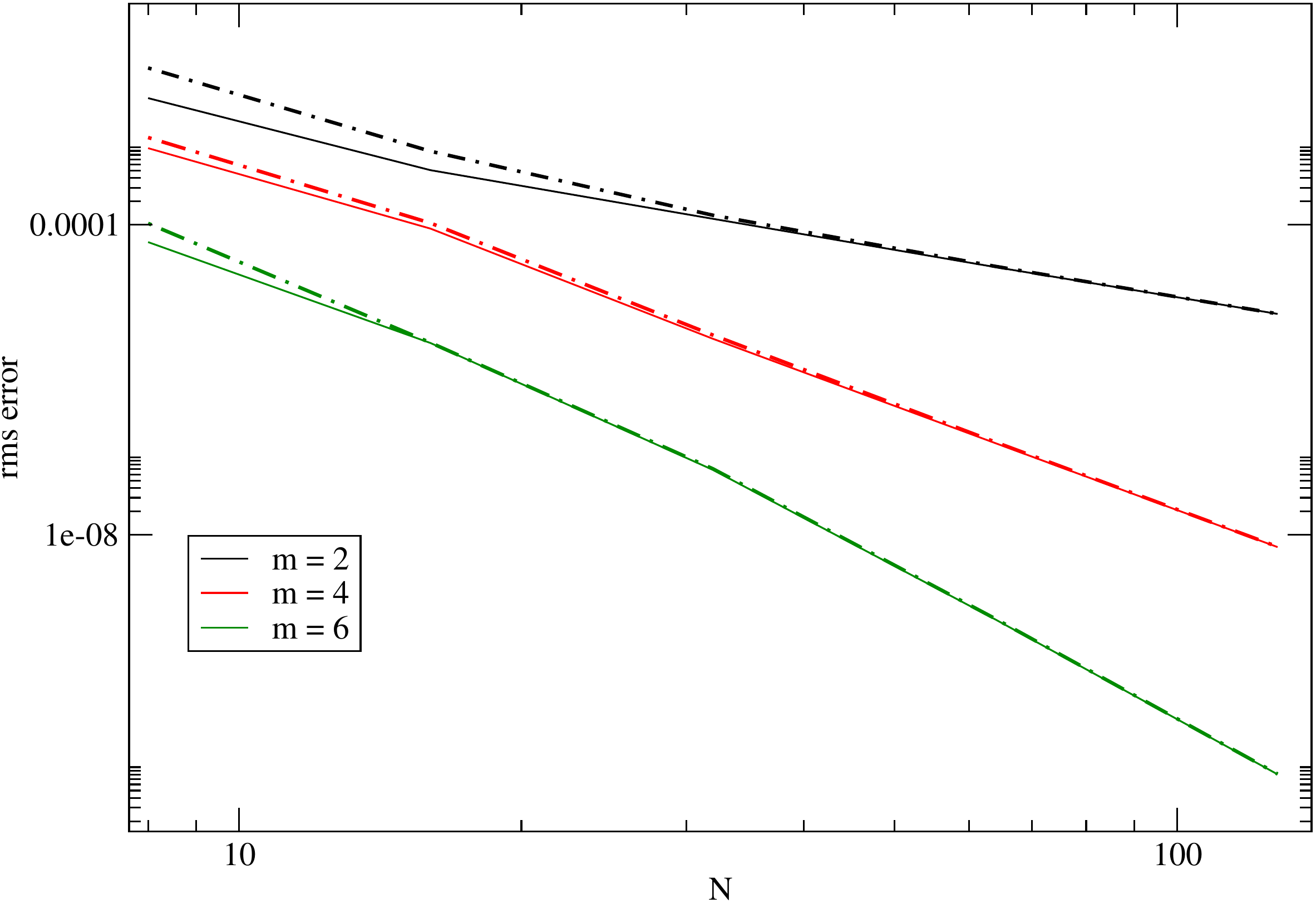}
  \caption{Convergence in discrete $l_2$-norm for the smooth manufactured solution in two-dimensions with Dirichlet boundary conditions using second, fourth and sixth order for collocated GMLS (dashed lines) and staggered GMLS (solid lines).}
  \label{manufacturedCompare}
\end{figure}
Figure \ref{manufacturedCompare} demonstrates the convergence in the $l_2$-norm of the two approaches for increasing refinement and increasing polynomial order. Both schemes provide $m^{th}$-order convergence, with the staggered scheme being slightly more accurate.

\begin{figure}[h!]
  \centering
  \includegraphics[width=0.45\textwidth]{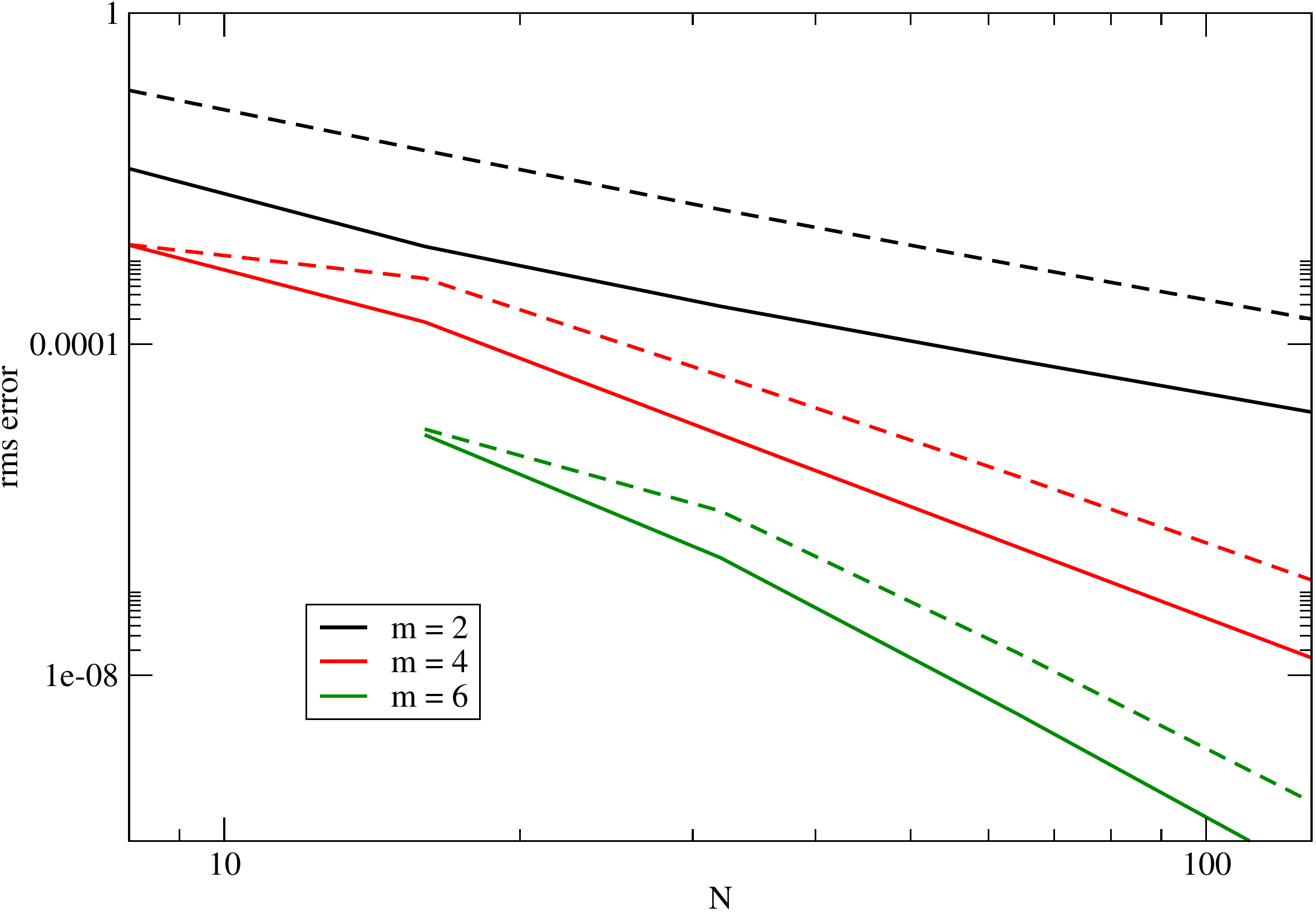}
  \includegraphics[width=0.45\textwidth]{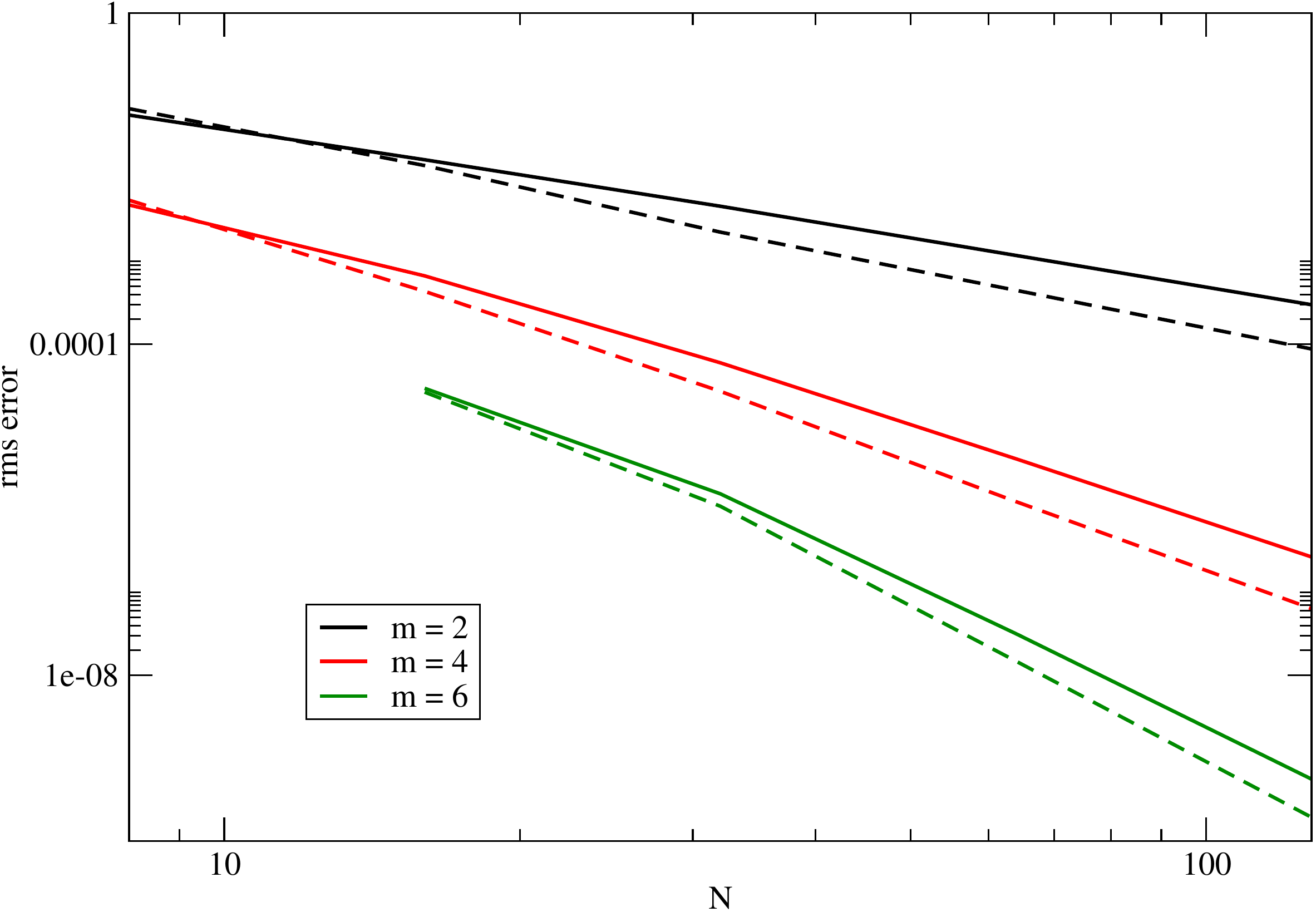}
  \caption{Convergence in discrete $l_2$-norm (left) and  $h_1$-seminorm (right) for manufactured solution with Dirichlet (solid lines) and Neumann (dashed lines) boundary conditions, suggesting a convergence result of $||u-u_{ex}||_{l_2} \leq C h^{m}$ and $||u-u_{ex}||_{h^1} \leq C h^{m}$  for the staggered scheme.}
  %
  %
\label{manufacturedL2H1}
\end{figure}

 Next, we examine the convergence rates of the staggered GMSL for both Dirichlet and Neumann boundary conditions. The data in Table \ref{tab:convrates} and the plots in Figure \ref{manufacturedL2H1} reveal equal convergence rates in both the $\ell_2$ norm and the $h^1$ seminorm, i.e., the scalar variable and its flux in the staggered GMLS converge at identical rates. This kind of behavior is typical of mixed methods for \eqref{eq:model} and lends further credence to the observations in Remark \ref{rem:compatible} about the similarities between our scheme and compatible discretization methods. Table \ref{tab:convrates3D} confirms that similar results hold in three dimensions.
%
\begin{table}[]
\centering
\caption{Convergence rate for the smooth manufactured solution problem in two-dimensions.}
\begin{tabular}{cccc}
\textit{\textbf{}}      & \textbf{m=2} & \textbf{m=4} & \textbf{m=6} \\ \hline
{$\ell_2$ - Dirichlet} & 2.085       & 4.470       &  6.486       \\
{$\ell_2$ - Neumann}   & 2.169       & 4.187       &  6.021       \\
{$h^1$ - Dirichlet}    & 1.979       & 3.940       &  5.839       \\
{$h^1$ - Neumann}      & 2.350       & 4.293       &  6.272      
\end{tabular}
\label{tab:convrates}
\end{table}

\begin{table}[]
\caption{Convergence rate for the smooth manufactured solution problem in three-dimensions and Dirichlet boundary conditions.}
\centering
\begin{tabular}{cccc}
\textit{\textbf{}}      & \textbf{m=2} & \textbf{m=4} & \textbf{m=6} \\ \hline
{$\ell_2$ } &  2.001       &  4.548       &  6.339       \\
{$h^1$ }    &  2.008       &  3.813       &  5.742       \\
\end{tabular}
\label{tab:convrates3D}
\end{table}

\subsection{Elliptic problems with rough coefficients}\label{sec:rough}
A hallmark of compatible discretization methods for \eqref{eq:model}, such as mixed finite elements, is their ability to correctly represent the flux variable for problems with discontinuous coefficients $\mu$. In such cases the normal component of the flux remains continuous across the material interface, while its tangential component is allowed to develop a discontinuity. Collocated methods, including Galerkin, stabilized Galerkin and least-squares finite elements cannot capture this behavior and tend to develop oscillations across the interface; see \cite{hughes} and \cite{bochevLSFEM} for further details and examples. 

Remark \ref{rem:compatible} points out structural similarities between the staggered GMLS and other compatible discretizations. The examples in this section aim to demonstrate that these similarities are not just formal  and that the new scheme is in fact capable of delivering physically correct approximations of the flux for problems with discontinuous coefficients. 

Our first example involves the so called five strip problem \cite{hughes}, which is a standard manufactured solution test case for examining the ability of a scheme to maintain normal flux continuity. In this problem $\phi_{ex} = 1-x$, $\Gamma_N=\Gamma$, the computational domain $\Omega=[0,1]^2$ is divided into five equal strips
\begin{equation}
\Omega_i = \left\{(x,y)\;|\; 0.2(i-1)\leq y \leq 0.2 i ; 0 \leq x \leq 1\right\},\quad i=1,\ldots,5
\end{equation}
and $\mu$ is assigned a different constant value $\mu_i$ on each $\Omega_i$. Here we use  $\mu_1 = 16$, $\mu_2 = 6$, $\mu_3 = 1$, $\mu_4 = 10$, and $\mu_5 = 2$. Substitution of $\phi_{ex}$ into \eqref{eq:model} defines the source term and the Neumann data $g$. In particular, since $\phi_{ex}$ is globally linear, $f=0$, $g=\pm \mu_i$ on the vertical parts of $\Gamma$, and
$$
\bm{u}\big|_{\Omega_i} = -\mu_i\nabla\phi_{ex} = \left(\begin{array}{c} \mu_i\\ 0\end{array}\right).
$$
Thus, the normal flux is continuous across the interface between two strips while the tangential component of $\bm{u}$ is a piecewise constant equal to $\mu_i$ on strop $\Omega_i$.

\begin{figure}[h!]
  \centering
  \includegraphics[width=0.95\textwidth]{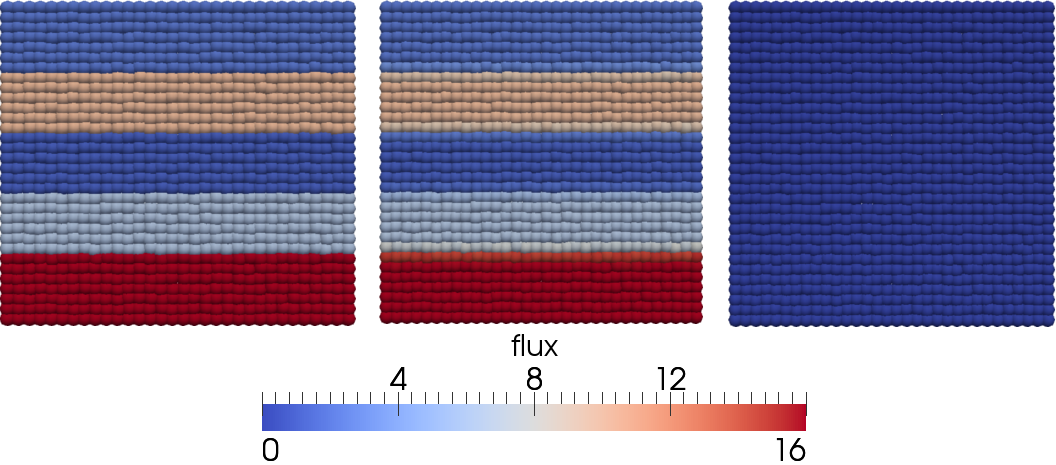}
  \caption{Fluxes for the Darcy flow strips case for non-conforming particle arrangement: exact longitudinal flux (left), staggered longitudinal flux (center), and staggered transverse flux (right).}
  \label{fluxplots}
\end{figure}

\begin{figure}[h!]
  \centering
  \includegraphics[width=0.45\textwidth]{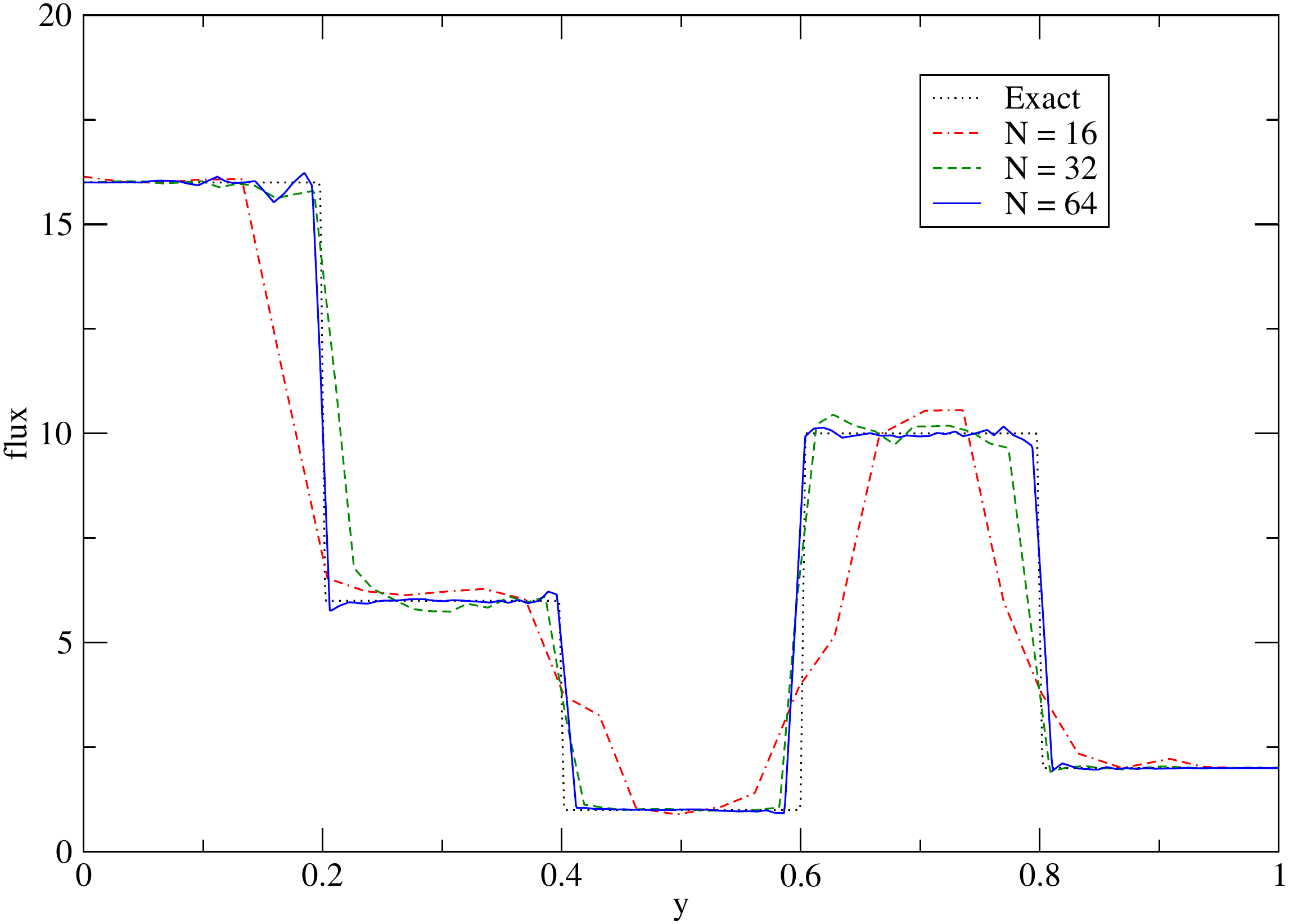}
  \includegraphics[width=0.45\textwidth]{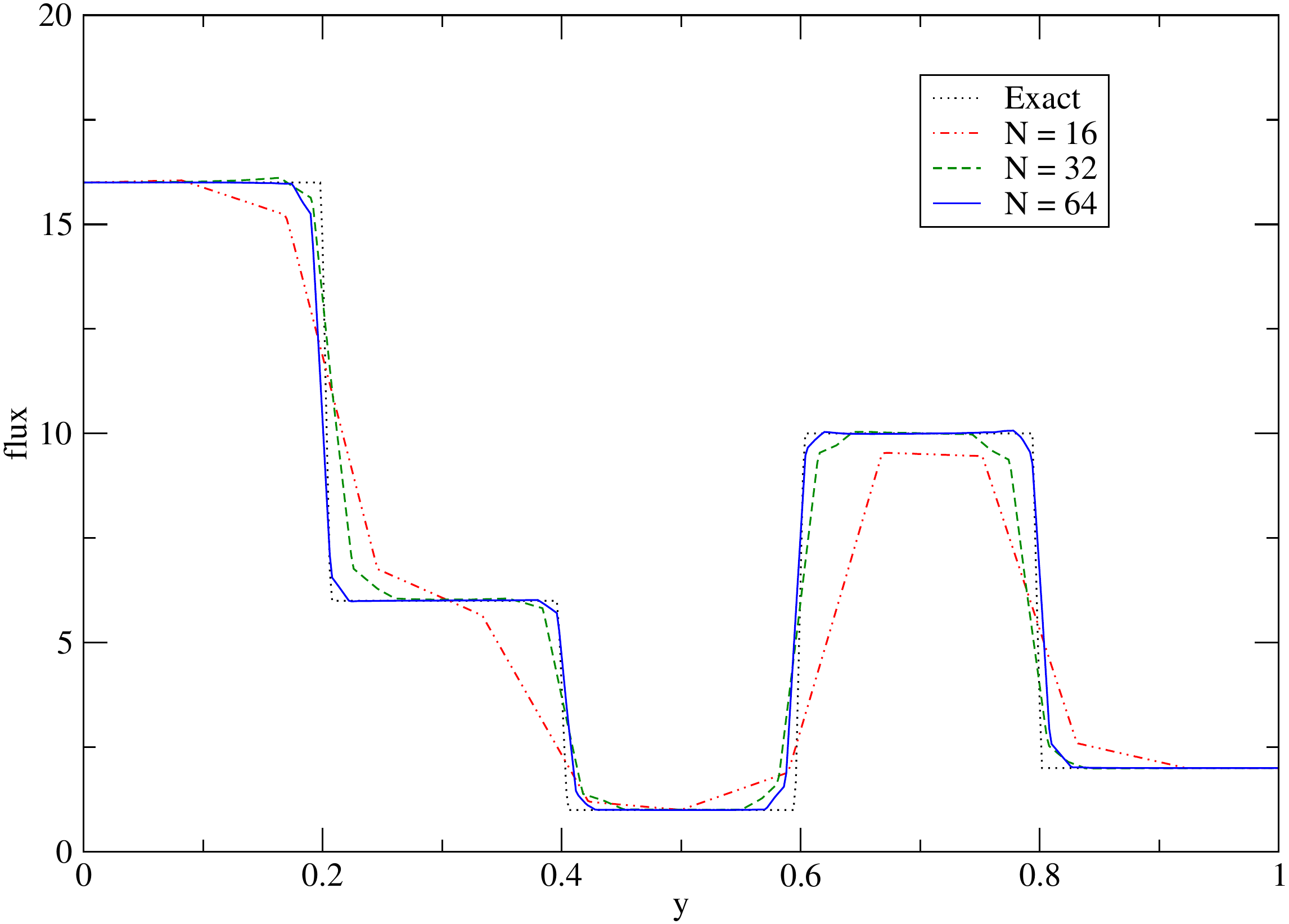}
  \caption{Transverse fluxes along the line $x = 0.5$ for the collocated (left) and staggered (right) schemes with increasing resolution.}
  \label{fluxplots2}
\end{figure}

The plots in Figure \ref{fluxplots} demonstrate that the staggered GMLS solution provides an accurate approximation of the tangential flux. Furthermore, comparison with the standard, collocated GMLS solution in Figure \ref{fluxplots2} reveals that the latter exhibits the same type of oscillations across the interfaces as found in nodal mesh-based methods \cite{hughes}.  It is worth pointing out that the staggered GMLS results in Figures \ref{fluxplots}--\ref{fluxplots2} were obtained without requiring that particles conform to the material interfaces between the strips. In contrast, mesh-based compatible methods typically require interface-fitted grids for accuracy. 

The model problem \eqref{eq:model} can also be used to study electrostatics problems in the absence of charge by taking $f=0$ and identifying $\mu$ as the electric permitivity (usually denoted $\epsilon$). In multiphysics applications, electromagnetic effects are often solved concurrently with a model for the mechanics governing domain deformation (e.g. electrophoresis of colloidal suspensions). When non-trivial geometry is considered for these applications, typically the only available tools to handle deforming domains with discontinuous material properties involve either costly arbitrary Lagrangian-Eulerian methods or some diffuse interfacial treatment \cite{SPMxianluo}. The current approach, on the other hand, allows a natural treatment of jumps in material properties without the need to develop a mesh conforming to discontinuities. 

\begin{figure}[h!]
  \centering
  \includegraphics[width=0.365\textwidth]{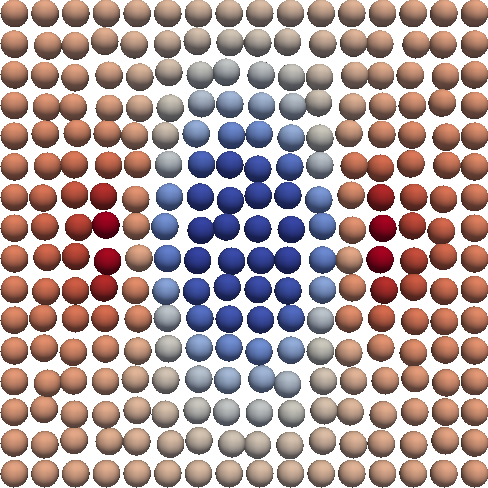}
  \includegraphics[width=0.45\textwidth]{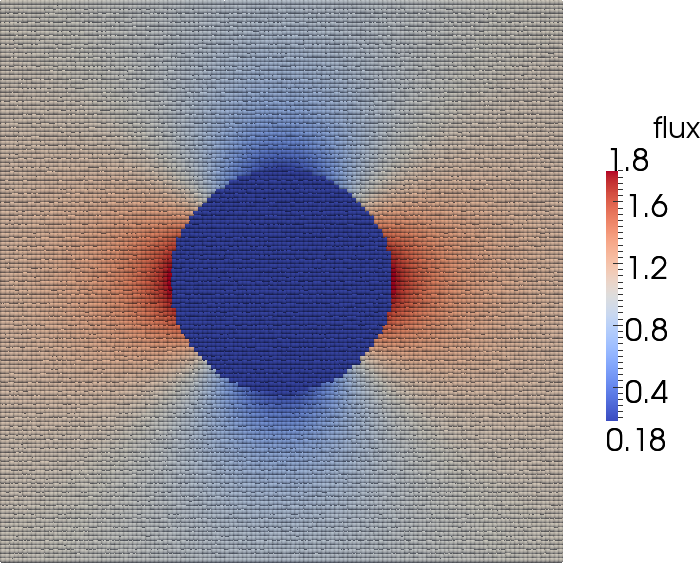}
  \caption{Plot of electric flux $\nabla \phi$ for coarsest ($N=16$) and finest ($N=256$) resolution.}
  \label{colloid}
\end{figure}

\begin{figure}[h!]
  \centering
  \includegraphics[width=0.45\textwidth]{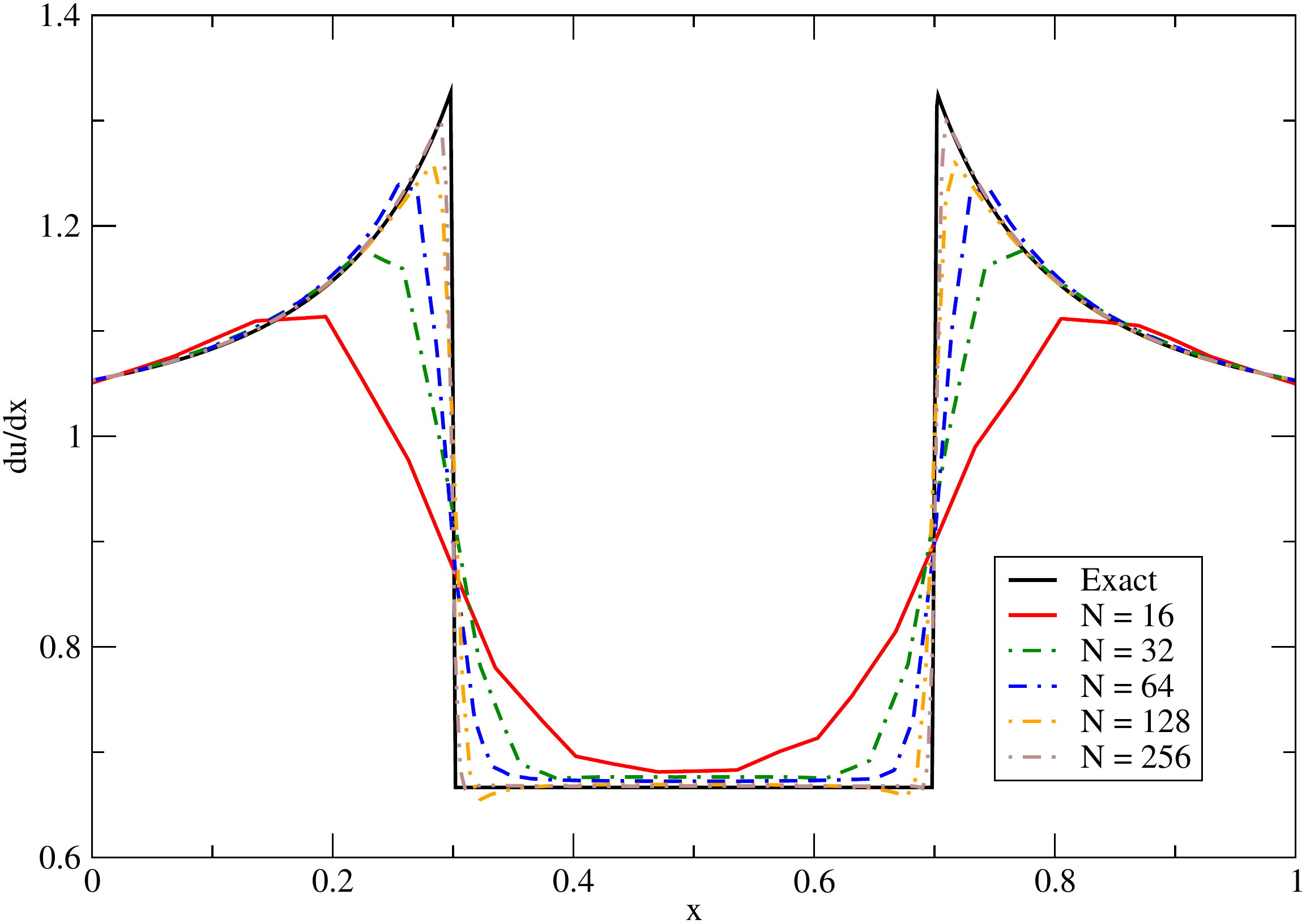}
  \includegraphics[width=0.45\textwidth]{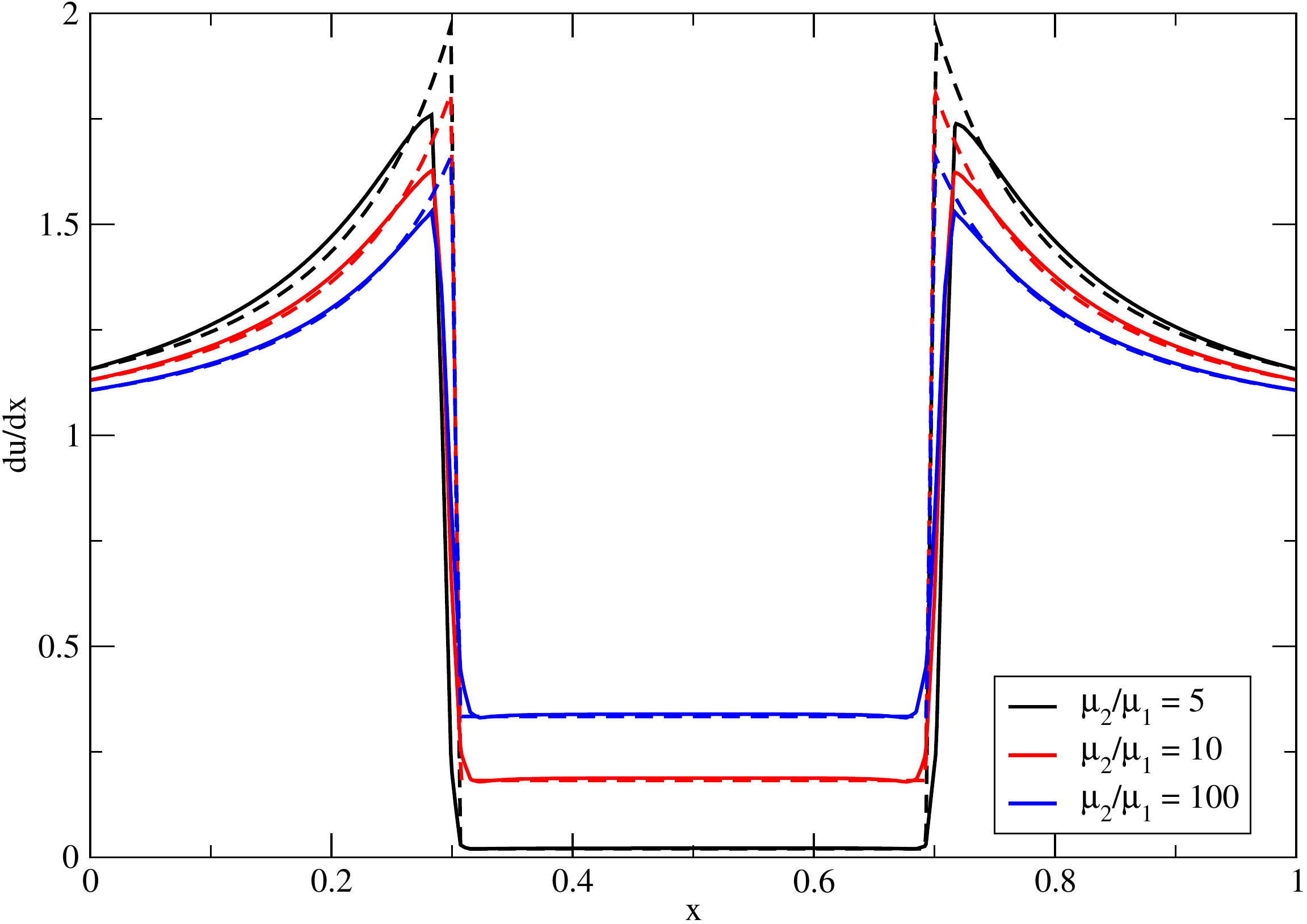}
  \caption{Electric flux along center line for increasing resolution (left) and for a range of permitivities (right).}
  \label{colloid2}
\end{figure}

Particles can be identified as belonging to either material with $\mu$ assigned appropriately. As an example, we simulate a circular cylinder of radius $R=1/2$ in a medium with permitivity $\mu/\mu_\infty = 2$ exposed to a uniform electric field $\nabla \phi = \left<1,0\right>$ for which an analytic solution can easily be calculated in an infinite domain. To solve numerically, we impose the analytic solution as a Dirichlet boundary condition on a unit square and demonstrate convergence in Figures \ref{colloid} and \ref{colloid2}. This discretization is stable for a ratio of permitivities ranging from $1-1000$ (Figure \ref{colloid2}). 


\begin{figure}[h!]
  \centering
  \includegraphics[width=0.45\textwidth]{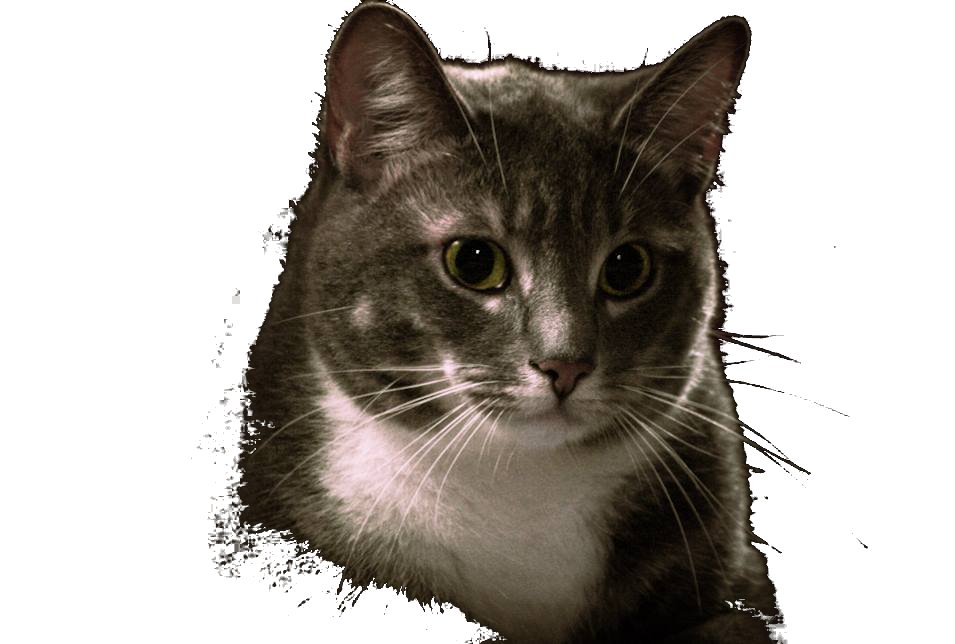}
  \includegraphics[width=0.45\textwidth]{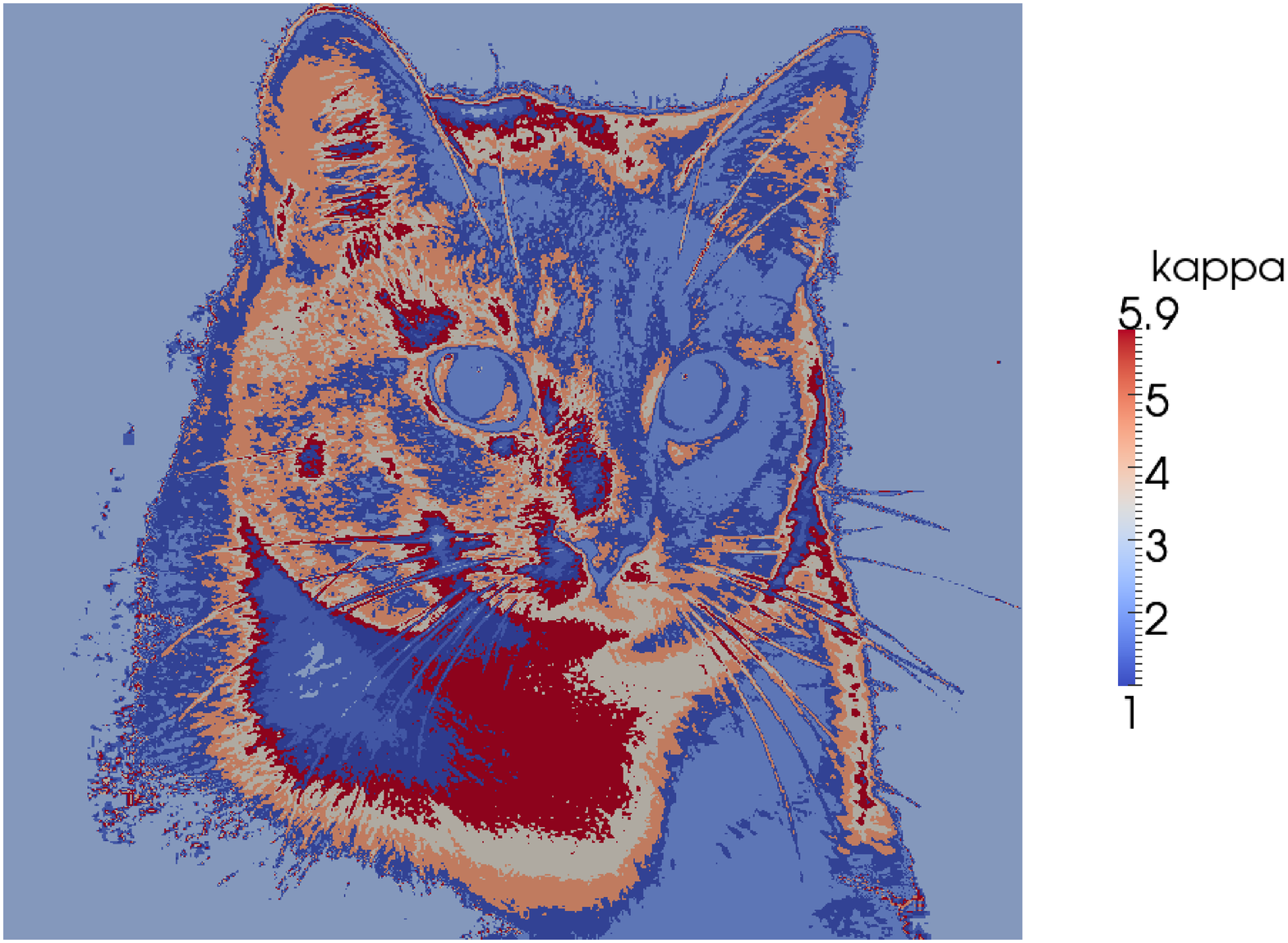}
  \caption{Conductivity for the underresolved stability test case.}
  \label{kappacat}
\end{figure}

As an example of the flexibility of the method for complex geometries, we now solve \eqref{eq:model} where $\mu$ is taken by processing an image of one of the authors' cat (Figure \ref{kappacat}) and assigning permitivities to regions of similar color. Particles are assigned values for $\mu_i$ by sampling pixels from the picture, and an approximately uniform electric field is imposed by setting Dirichlet boundary conditions.
\begin{equation}
\phi_i|_{\partial \Omega} = 1 - x_i
\end{equation}

\begin{figure}[h!]
  \centering
  \includegraphics[width=1.0\textwidth]{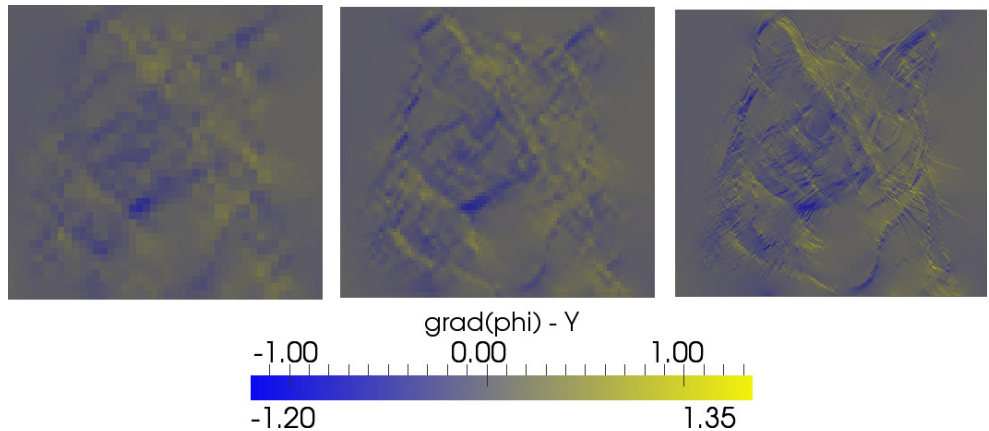}
  \caption{$\partial_y \phi$ for underresolved stability test case. Underresolved cases provided similar results despite highly oscillatory conductivity below resolution lengthscale.}
  \label{refinedcat}
\end{figure}

The robustness of the method in the presence of fluctuations in the permitivity below the resolution lengthscale is shown in the refinement study presented in Figure \ref{refinedcat}. Examining the y-component of the resulting ``electric field'' for resolution ranging from $32^2$ to $1024^2$, the results remain consistent and non-oscillatory. We note that, using AMG, the resulting matrix solve took only a few seconds on a desktop computer. While this is a slightly whimsical demonstration of the flexibility of the method, we hope to motivate that this approach could easily be adopted to any problem in which voxel data is available of the underlying diffusivity (e.g. applications in MRI imaging, geological imaging, etc).


\section{Conclusions and future work}\label{sec:conclude}
Motivated by mimetic methods on primal-dual grids we have developed a meshfree discretization for a model diffusion equation, which exhibits many of the attractive computational properties of compatible discretization. Our scheme generalizes staggered finite difference methods to arbitrary point clouds and particle arrangements by using separate discretizations of the divergence and gradient operators defined locally at the neighborhood of each discretization point.
While traditional staggered methods are of either finite difference type and are restricted to Cartesian grids, or are of finite volume type and require a primal-dual grid complex, this approach is truly meshfree since it  uses only 
the $\epsilon$-neighborhood graph of particle connectivity. 

Due to the lack of a mesh, there is no obvious inner-product space in which to analyze this method and we have opted instead to demonstrate the compatible nature of the scheme by comparison to a series of problems in which discretizations lacking compatibility are prone to failure.

For the sake of brevity, we have restricted the focus of this work to the model problem \eqref{eq:model}. We note here that the stability and high-order accuracy achieved for this problem carries over to the Stokes problem, which we will discuss in a later work.
\section{Acknowledgments}
{\it This material is based upon work supported by the U.S. Department of Energy Office of Science, Office of Advanced Scientific Computing Research, Applied Mathematics program as part of the Colloboratory on Mathematics for Mesoscopic Modeling of Materials (CM4), under Award Number DE-SC0009247. This research used resources of the National Energy Research Scientific Computing Center, which is supported by the Office of Science of the U.S. Department of Energy under Contract No. DE-AC02-05CH11231. Sandia is a multiprogram laboratory operated by Sandia Corporation, a Lockheed Martin Company, for the U.S. Department of Energy under contract DE-AC04-94-AL85000.
}

\bibliographystyle{ieeetr}
\bibliography{refs}

\end{document}